\documentclass[11pt]{article}

\usepackage{amsthm}
\usepackage{amsmath}
\usepackage{amssymb}
\usepackage{graphicx}

\newtheorem{theorem}{Theorem}[section]
\newtheorem{lemma}[theorem]{Lemma}
\newtheorem{question}[theorem]{Question}
\newtheorem{conjecture}[theorem]{Conjecture}

\newtheorem{proposition}[theorem]{Proposition}
\newtheorem{corollary}[theorem]{Corollary}

\newtheorem{problem}[theorem]{Problem}

\newcommand{\A}{{\cal A}}
\newcommand{\B}{{\cal B}}
\newcommand{\C}{\mathbb{C}}
\newcommand{\E}{E_{\bullet}}
\newcommand{\F}{F_{\bullet}}
\newcommand{\G}{G_{\bullet}}
\renewcommand{\H}{\mathcal{H}}
\newcommand{\I}{\mathcal{I}}
\renewcommand{\L}{\mathcal{L}}
\newcommand{\R}{\mathbb{R}}

\newcommand{\rk}{\textrm{rk}}
\newcommand{\s}{\textrm{size}}
\newcommand{\T}{\mathcal{T}}
\newcommand{\Fl}{\mathcal{F}\ell}
\newcommand{\Z}{\mathbb{Z}}
\newcommand{\cinum}[1]{\put(5,5){\circle{15}} #1}

\newlength{\cellsize} 
\newsavebox{\cell}
\sbox{\cell}{\begin{picture}(18,18) \put(0,0){\line(1,0){18}}
\put(0,0){\line(0,1){18}} \put(18,0){\line(0,1){18}}
\put(0,18){\line(1,0){18}}
\end{picture}}
\newcommand\cellify[1]{\def\thearg{#1}\def\nothing{}%
\ifx\thearg\nothing \vrule width0pt height\cellsize depth0pt\else
\hbox to 0pt{\usebox{\cell} \hss}\fi%
\vbox to \cellsize{ \vss \hbox to \cellsize{\hss$#1$\hss} \vss}}
\newcommand\tableau[1]{\vtop{\let\\\cr
\baselineskip -16000pt \lineskiplimit 16000pt \lineskip 0pt
\ialign{&\cellify{##}\cr#1\crcr}}}
\newcommand{\emt}{\mbox{ }}

\def\blfootnote{\xdef\@thefnmark{}}

\begin{document}

\title{Flag arrangements and \\ triangulations of products of simplices.}
\author{Federico Ardila \and Sara Billey\footnote{Both
authors supported by NSF grant DMS-9983797. \newline \textit{Date:}
\today \newline \textit{Keywords:} matroids, permutation
arrays, Schubert calculus, tropical hyperplane arrangements,
Littlewood-Richardson coefficients}}
\date{}
\maketitle

\begin{abstract}
We investigate the line arrangement that results from intersecting
$d$ complete flags in $\C^n$. We give a combinatorial description
of the matroid $\T_{n,d}$ that keeps track of the linear
dependence relations among these lines.

We prove that the bases of the matroid $\T_{n,3}$ characterize the
triangles with holes which can be tiled with unit rhombi. More
generally, we provide evidence for a conjectural connection
between the matroid $\T_{n,d}$, the triangulations of the product
of simplices $\Delta_{n-1} \times \Delta_{d-1}$, and the
arrangements of $d$ tropical hyperplanes in tropical
$(n-1)$-space.

Our work provides a simple and effective criterion to ensure the
vanishing of many Schubert structure constants in the flag
manifold, and a new perspective on Billey and Vakil's method for
computing the non-vanishing ones.
\end{abstract}

\section{Introduction.}
\label{section:intersecting.flags}

Let $\E^1, \ldots, \E^d$ be $d$ generically chosen complete flags
in $\C^n$. Write
\[
\E^k = \{\{0\}=E^k_0 \subset E^k_1 \subset \cdots \subset E^k_n =
\C^n\},
\]
where $E^k_i$ is a vector space of dimension $i$. Consider the set
${\bf E}_{n,d}$ of one-dimensional intersections determined by the
flags; that is, all lines of the form $E^1_{a_1} \cap E^2_{a_2}
\cap \cdots \cap E^d_{a_d}$.

The initial goal of this paper is to characterize the line
arrangements $\C^n$ which arise in this way from $d$ generically
chosen complete flags. We will then show an unexpected connection
between these line arrangements and an important and ubiquitous
family of subdivisions of polytopes: the triangulations of the
product of simplices $\Delta_{n-1} \times \Delta_{d-1}$. These
triangulations appear naturally in studying the geometry of the
product of all minors of a matrix \cite{Babson}, tropical geometry
\cite{Develin}, and transportation problems \cite{Sturmfels}. To
finish, we will illustrate some of the consequences that the
combinatorics of these line arrangements have on the Schubert
calculus of the flag manifold.

\medskip

The results of the paper are roughly divided into four parts as
follows. First of all, Section \ref{section:H_mn} is devoted to
studying the line arrangement determined by the intersections of a
generic arrangement of hyperplanes. This will serve as a warmup
before we investigate generic arrangements of complete flags, and
the results we obtain will be useful in that investigation.

The second part consists of Sections
\ref{section:flags.to.simplex}, \ref{section:T_nd} and
\ref{section:T_nd.is.right}, where we will characterize the line
arrangements that arise as intersections of a ``matroid-generic"
arrangement of $d$ flags in $\C^n$. Section
\ref{section:flags.to.simplex} is a short discussion of the
combinatorial setup that we will use to encode these geometric
objects. In Section \ref{section:T_nd}, we propose a combinatorial
definition of a matroid $\T_{n,d}$. In Section
\ref{section:T_nd.is.right} we will show that $\T_{n,d}$ is the
matroid of the line arrangement of any $d$ flags in $\C^n$ which
are generic enough. Finally, we show that these line arrangements
are completely characterized combinatorially: any line arrangement
in $\C^n$ whose matroid is $\T_{n,d}$ arises as an intersection of
$d$ flags.

The third part establishes a surprising connection between these
line arrangements and an important class of subdivisions of
polytopes. The bases of $\T_{n,3}$ exactly describe the ways of
punching $n$ triangular holes into the equilateral triangle of size
$n$, so that the resulting holey triangle can be tiled with unit
rhombi. A consequence of this is a very explicit geometric
representation of $\T_{n,3}$. We show these results in Section
\ref{section:rhombus.tilings}. We then pursue a higher-dimensional
generalization of this result. In Section
\ref{section:fine.mixed.subdivisions}, we suggest that the fine
mixed subdivisions of the Minkowski sum $n \Delta_{d-1}$ are an
adequate $(d-1)$-dimensional generalization of the rhombus tilings
of holey triangles. We give a completely combinatorial description
of these subdivisions. Finally, in Section
\ref{section:triangulations.matroid}, we prove that each pure mixed
subdivision of the Minkowski sum $n\Delta_{d-1}$ (or equivalently,
each triangulation of the product of simplices $\Delta_{n-1} \times
\Delta_{d-1}$) gives rise to a basis of $\T_{n,d}$. We conjecture
that every basis of $\T_{n,d}$ arises in this way. In fact, we
conjecture that every basis of $\T_{n,d}$ arises from a
\emph{regular} subdivision or, equivalently, from an arrangement of
$d$ tropical hyperplanes in tropical $(n-1)$-space.

The fourth and last part of the paper, Section
\ref{section:schubert}, presents some of the consequences of our
work in the Schubert calculus of the flag manifold. We start by
recalling Eriksson and Linusson's permutation arrays, and Billey
and Vakil's related method for explicitly intersecting Schubert
varieties. In Section \ref{subsection:genericity} we show how the
geometric representation of the matroid $\T_{n,3}$ of Section
\ref{section:rhombus.tilings} gives us a new perspective on Billey
and Vakil's method for computing the structure constants $c_{uvw}$
of the cohomology ring of the flag variety. Finally, Section
\ref{subsection:schubert.zero} presents a simple and effective
criterion for guaranteeing that many Schubert structure constants
are equal to zero.

\medskip

We conclude with some future directions of research that are
suggested by this project. %In particular, there is still much work
%to be done in clarifying the relationship between the two main
%subjects of our paper: the geometry of $d$ flags in $\C^n$ and the
%triangulations of $\Delta_{n-1} \times \Delta_{d-1}$. In this
%paper, we show that some aspects of the geometric information of
%the flags (the combinatorics of the line arrangement they
%determine) are described in a small set of tiles of the
%triangulations (the $n$ ``pure" tiles). It would be extremely
%interesting to extend this connection, and explain how more subtle
%geometric questions about flags relate to the complete
%triangulations, and to their multiple appearances in tropical
%geometry, optimization, and other subjects.

\section{The lines in a generic hyperplane arrangement.} \label{section:H_mn}

Before thinking about flags, let us start by studying the slightly
easier problem of understanding the matroid of lines of a generic
arrangement of $m$ hyperplanes in $\C^n$. We will start by
presenting, in Proposition \ref{prop:H_mn}, a combinatorial
definition of this matroid $\H_{n,m}$. Theorem
\ref{thm:generichyperplanes} then shows that this is, indeed, the
right matroid. As it turns out, this warmup exercise will play an
important role in Section \ref{section:T_nd.is.right}.

Throughout this section, we will consider an arrangement of $m$
generically chosen hyperplanes $H_1, \ldots, H_m$ in $\C^n$ passing
through the origin. For each subset $A$ of $[m]=\{1,2,\ldots, m \}$,
let
\[
H_A = \bigcap_{a \in A} H_a.
\]
By genericity,
\begin{displaymath}
\dim H_A = \left\{ \begin{array}{ll} n-|A| & \textrm{if $|A| \leq
n$,} \\
0 & \textrm{otherwise.}
\end{array} \right.
\end{displaymath}
Therefore, the set $L_{n,m}$ of one-dimensional intersections of
the $H_i$s consists of the ${m \choose n-1}$ lines $H_A$ for
$|A|=n-1$.

There are several ``combinatorial" dependence relations among the
lines in $L_{n,m}$, as follows. Each $t$-dimensional intersection
$H_B$ (where $B$ is an $(n-t)$-subset of $[m]$) contains the lines
$H_A$ with $B \subseteq A$. Therefore, in an independent set
$H_{A_1}, \ldots, H_{A_k}$ of $L_{n,m}$, we cannot have $t+1$
$A_i$s which contain a fixed $(n-t)$-set $B$.

At first sight, it seems intuitively clear that, in a generic
hyperplane arrangement, these will be the only dependence
relations among the lines in $L_{n,m}$. This is not as obvious as
it may seem: let us illustrate a situation in $L_{4,5}$ which is
surprisingly close to a counterexample to this statement. For
simplicity, we will draw the three-dimensional projective picture.
Each hyperplane in $\C^4$ will now look two-dimensional, and the
lines in the arrangement $L_{4,5}$ will look like points. Denote
hyperplanes $H_1, \ldots, H_5$ simply by $1,\ldots,5$, and an
intersection like $H_{124}$ simply by $124$.

\begin{figure}[h]
\centering
\includegraphics[height=5cm]{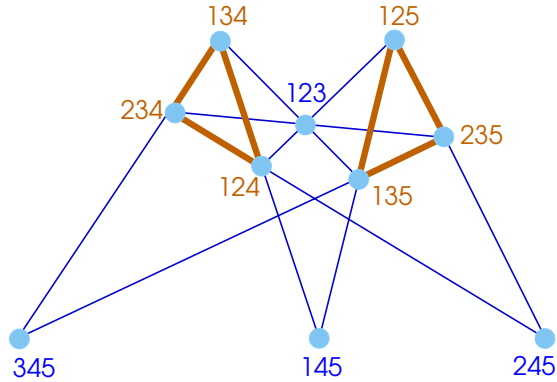}
\caption{The Desargues configuration in $L_{4,5}$.}
\label{fig:desargues}
\end{figure}

In Figure \ref{fig:desargues}, we have started by drawing the
triangles $T$ and $T'$ with vertices $124,234,134$ and
$125,235,135$, respectively. The three lines connecting the pairs
$(124,125)$, $(234,235)$ and $(134,135)$, are the lines $12, 23,$
and $13$, respectively. They intersect at the point $123$, so that
the triangles $T$ and $T'$ are perspective with respect to this
point.

Now Desargues' theorem applies, and it predicts an unexpected
dependence relation. It tells us that the three points of
intersection of the corresponding sides of $T$ and $T'$ are
collinear. The lines $14$ (which connects $124$ and $134$) and
$15$ (which connects $125$ and $135$) intersect at the point
$145$. Similarly, $24$ and $25$ intersect at $245$, and $34$ and
$35$ intersect at $345$. Desargues' theorem says that the points
$145$, $245$, and $345$ are collinear. In principle, this new
dependence relation does not seem to be one of our predicted
``combinatorial relations". Somewhat surprisingly, it is: it
simply states that these three points are on the line $45$.

The previous discussion illustrates two points. First, it shows
that Desargues' theorem is really a combinatorial statement about
incidence structures, rather than a geometric statement about
points on the Euclidean plane. Second, and more important to us,
it shows that even five generic hyperplanes in $\C^4$ give rise to
interesting geometric configurations.
%the very simple arrangement of lines $L_{5,4}$ already contains
%configurations of points which may have non-trivial geometric
%properties.
It is not unreasonable to think that larger arrangements $L_{n,m}$
will contain other configurations, such as the Pappus
configuration, which have nontrivial and honestly geometric
dependence relations that we may not have predicted.

Having told our readers what they might need to worry about, we
now intend to convince them not to worry about it.

%Having seen what problems we might need to worry about, we will
%now show that those problems do not arise in the situation at
%hand.

First we show that the combinatorial dependence relations in
$L_{n,m}$ are consistent, in the sense that they define a matroid.

\begin{proposition}\label{prop:H_mn}
Let $\I$ consist of the collections $I$ of subsets of $[m]$, each
containing $n-1$ elements, such that no $t+1$ of the sets in $I$
contain an $(n-t)$-set. In symbols,
%\[
%\I := \left\{\{I_1, \ldots, I_k\} \subseteq {[m] \choose n-1} |
%\textrm{ for } \{i_1, \ldots, i_s\} \subset [k], |I_{i_1} \cap
%\cdots \cap I_{i_{s}}| \leq n-s\right\}.
%\]%
%\[
%\I := \left\{\{I_1, \ldots, I_k\} \subseteq {[m] \choose n-1}
%\large{|} \textrm{ for all } S \subseteq [k], |\bigcap_{s \in S}
%I_{s}| \leq n-|S|\right\}.
%\]
\[
\I := \left\{I \subseteq {[m] \choose n-1} \textrm{ such that for
all } S \subseteq I, \,\,\, \bigl|\hspace{-.1cm}\bigcap_{A \in S}
A\bigr| \, \leq n-|S|\right\}.
\]
Then $\I$ is the collection of independent sets of a matroid
$\H_{n,m}$.
\end{proposition}

\begin{proof}
A circuit of that matroid would be a minimal collection $C$ of $s$
subsets of $[m]$ of size $n-1$, all of which contain one fixed
$(n-s+1)$-set. It suffices to verify the circuit axioms:

\noindent{\bf (C1)} No proper subset of a circuit is a circuit.

\noindent{\bf (C2)} If two circuits $C_1$ and $C_2$ have an
element $x$ in common, then $C_1\cup C_2 -x$ contains a circuit.

The first axiom is satisfied trivially. Now consider two circuits
$C_1$ and $C_2$ containing a common $(n-1)$-set $X_1$. Let
\[
C_1 = \{X_1, \ldots, X_a, Y_1, \ldots, Y_b\}, \qquad C_2 = \{X_1,
\ldots, X_a, Z_1, \ldots, Z_c\},
\]
where the $Y_i$s and $Z_i$s are all distinct. Write
\[
X = \bigcap_{i=1}^a X_i,  \qquad Y = \bigcap_{i=1}^b Y_i,  \qquad
Z = \bigcap_{i=1}^c Z_i.
\]
By definition of $C_1$ and $C_2$ we have that $|X \cap Y| \geq
n-(a+b)+1$ and $|X \cap Z| \geq n-(a+c)+1$, and their minimality
implies that $|X| \leq n-a$. Therefore
\begin{eqnarray*}
|X \cap Y \cap Z| &=& |X \cap Y| + |X \cap Z| - |(X \cap Y) \cup
(X \cap Z)| \\
&\geq& |X \cap Y| + |X \cap Z| - |X| \\
&\geq& (n-a-b+1) + (n-a-c+1) - (n-a) \\
&=& n-a-b-c+2,
\end{eqnarray*}
and hence
\[
|X_2 \cap \cdots \cap X_a \cap Y_1 \cap \cdots \cap Y_b \cap Z_1
\cap \cdots \cap Z_c| \geq n - (a+b+c-1) + 1.
\]
It follows that $C_1 \cup C_2 - X_1$ contains a circuit, as
desired.

%
%It is possible to check directly that $\I$ satisfies the axioms
%for a collection of independent sets of a matroid.
%
%It is perhaps easier to recall the following result. Given a set
%$E$ and a function $f:2^E \rightarrow \R$ which is increasing
%($f(X) \leq f(Y)$ for $X \subseteq Y$) and submodular ($f(X)+f(Y)
%\geq f(X \cup Y) + f(X \cap Y)$), there exists a matroid $M(f)$ on
%$E$ whose circuits are the minimal sets $S \subseteq E$ such that
%$|S| > f(S)$. This matroid is denoted $M(f)$, and called the
%\emph{polymatroid} of the function $f$. \cite[find]{Oxley}
%\textsf{Find reference, check details}.
%
%It is then immediate that $\H_{n,m}$ is the polymatroid
%corresponding to the function
%\[
%f(S) := n - |\hspace{-.1cm}\bigcap_{A \in S} A|,
%\]
%which is clearly increasing and submodular.
\end{proof}

Now we show that this matroid $\H_{n,m}$ is the one determined by
the lines in a generic hyperplane arrangement.

\begin{theorem}\label{thm:generichyperplanes}
If a central\footnote{A hyperplane arrangement is \emph{central}
if all its hyperplanes go through the origin.} hyperplane
arrangement $\A=\{H_1, \ldots, H_m\}$ in $\C^n$ is generic enough,
then the matroid of the ${m \choose n-1}$ lines $H_A$ is
isomorphic to $\H_{n,m}$.
\end{theorem}

\begin{proof}
We already observed that the one-dimensional intersections of $\A$
satisfy all the dependence relations of $\H_{n,m}$. Now we wish to
show that, if $\A$ is ``generic enough", these are the only
relations.

Construct a hyperplane arrangement as follows. Consider the $m$
coordinate hyperplanes in $\C^m$, numbered $J_1,\ldots, J_m$. Pick
a sufficiently generic $n$-dimensional subspace $V$ of $\C^m$, and
consider the ($(n-1)$-dimensional) hyperplanes $H_1 = J_1 \cap V,
\ldots, H_m = J_m \cap V$ in $V$. We will see that, if $V$ is
generic enough in the sense of Dilworth truncations, then the
arrangement $\{H_1, \ldots, H_m\}$ is generic enough for our
purposes. We now recall this setup.
%
%Consider a matroid $M$ represented as a collection of lines in
%$\C^r$, and consider a subspace $V$ of $\C^r$ of codimension $k-1$
%which is ``generic with respect to $M$". Each $k$-dimensional flat
%$F$ of $M$ intersects $V$ in a line $F \cap V$. When $V$ is
%``generic enough", this collection of lines gives rise to a fixed
%matroid $D_k(M)$, called the $k$-th Dilworth truncation of $M$.

\begin{theorem}(Brylawski, Dilworth, Mason, \cite{Brylawski, Brylawski2,
Mason})\label{thm:dilworth} Let $L$ be a set of lines in $\C^r$
whose corresponding matroid is $M$. Let $V$ be a subspace of
$\C^r$ of codimension $k-1$. For each $k$-flat $F$ spanned by $L$,
let $v_F = F \cap V$.
\begin{enumerate}
\item If $V$ is generic enough, then each $v_F$ is a line, and the
matroid $D_k(M)$ of the lines $v_F$ does not depend on $V$. \item
The circuits of $D_k(M)$ are the minimal sets $\{v_{F_1}, \ldots,
v_{F_a}\}$ such that $\rk_M(F_1 \cup \cdots \cup F_a) \leq
a+k-2$.\footnote{The idea behind this is that, if the span of
$F_1, \ldots, F_a$ has dimension less than $a+k-1$, then, upon
intersection with $V$ (which has codimension $k-1$), their span
will have dimension less than $a$.} This matroid is called the
\emph{$k$-th Dilworth truncation} of $M$.\footnote{The matroid
$D_k(M)$ can be defined combinatorially by specifying its circuits
in the same way, even if $M$ is not representable. In fact, when
$M$ is representable, the most subtle aspect of our definition of
$D_k(M)$ is the construction of a ``generic enough" subspace $V$,
and hence of a geometric realization of $D_k(M)$. This
construction was proposed by Mason \cite{Mason} and proved correct
by Brylawski \cite{Brylawski}. They also showed that, if $M$ is
not realizable, then $D_k(M)$ is not realizable either.}
\end{enumerate}
\end{theorem}
This is precisely the setup that we need. Let $L = \{1, \ldots,
m\}$ be the coordinate axes of $\C^m$, labelled so that coordinate
hyperplane $J_i$ is normal to axis $i$. These $m$ lines are a
realization of the free matroid $M_m$ on $m$ elements. %The
%$(m-1)$-flats of $M_m$ are the sets $[m]-i$; geometrically, they
%correspond to the coordinate hyperplanes $J_1, \ldots, J_m$. Each
%$(m-k)$-flat is combinatorially an $(m-k)$-subset $T$ of $[m]$,
%and geometrically an intersection of the $k$ hyperplanes $J_i$
%with $i \notin T$.

Now consider the $(m-n+1)$-th Dilworth truncation $D_{m-n+1}(M_m)$
of $M_m$, obtained by intersecting our configuration with an
$n$-dimensional subspace $V$ of $\C^m$, which is generic enough
for Theorem \ref{thm:dilworth} to apply. For each $(m-n+1)$-subset
$T$ of $L=\{1,\ldots,m\}$, we get an element of the matroid of the
form
%Combinatorially, each element of the matroid is just an
%$(m-n+1)$-subset $T$ of $[m]$. Geometrically, this element is
%represented by the line
\[
v_T = (\textrm{span } T) \cap V = \Bigl( \bigcap_{i \notin T} J_i
\Bigr) \cap V = \bigcap_{i \notin T} (J_i \cap V) = \bigcap_{i
\notin T} H_i = H_{[m]-T},
\]
where, as before, $H_i = J_i \cap V$ is a hyperplane in $V$. Since
$\bigl|[m]-T\bigr| = n-1$, this $v_T$ is precisely one of the
lines in the arrangement $L_{n,m}$ of one-dimensional
intersections of $\{H_1, \ldots, H_m\}$. In Theorem
\ref{thm:dilworth}, we have a combinatorial description for the
matroid $D_{m-n+1}(M_m)$ of the $v_T$s. It remains to check that
this matches our description of $\H_{n,m}$.

This verification is straightforward.
%This is precisely the setup that we need. Consider the $m$
%coordinate axes $x_1, \ldots, x_m$ of $\C^m$, which are a
%realization of the free matroid on $m$ elements. We considered the
%intersections of the $(m-1)-flats$ of this matroid (which are the
%coordinate hyperplanes $J_1, \ldots, J_m$) with a generic subspace
%$V$ of codimension $m-n+1$. Our aim was to study the matroid
%determined by the intersections
%\[
%H_I = \bigcap_{i \in I} H_i = \bigcap_{i \in I} (J_i \cap V),
%\]
%where $I$ is an $(n-1)$-subset of $[m]$. But
%\[
%H_I = \left(\bigcap_{i \in I} J_I\right)\cap V =
%\left(\textrm{span}\{x_i \, | \, i \notin I\}\right) \cap V,
%\]
%which is one of the elements of (this realization of) the
%$(m-n+1)$-th Dilworth truncation of the free matroid on $m$
%elements.
%
%So it suffices to check that the genericity that Dilworth
%truncation prescribes on the $H_S$s is the same as the genericity
%we claim to have.
In $D_{m-n+1}(M_m)$, the collection $\{v_{T_1}, \ldots, v_{T_a}\}$
is a circuit if it is a minimal set such that the following
equivalent conditions hold:
\begin{eqnarray*}
\rk_{M_m}(T_1 \cup \cdots \cup T_a) & \leq & a + (m-n+1) - 2, \\
|T_1 \cup \cdots \cup T_a| & \leq & m-(n-a+1), \\
|([m]-T_1) \cap \cdots \cap ([m]-T_a)| & \geq & n-a+1.
\end{eqnarray*}
This is equivalent to $\{[m]-T_1, \ldots, [m]-T_a\}$ being a
circuit of the matroid $\H_{n,m}$, which is precisely what we
wanted to show. This completes the proof of Theorem
\ref{thm:generichyperplanes}.
%
%
%In the Dilworth truncation, $\{H_{[m]-T_1}, \ldots, H_{[m]-T_a}\}$
%is a circuit if it is a minimal set such that
%\[ \rk\left( (\textrm{span}\{x_i \, | \, i \notin
%S_1\}) \cup \cdots \cup (\textrm{span}\{x_i \, | \, i \notin
%S_a\})\right)
% \leq a + (m-n+1) - 2,
% \]
%which, in the free matroid, is equivalent to
%\[
%|([m]-S_1) \cup \cdots \cup ([m]-S_a)| \leq m-(n-a+1),
%\]
%or
%\[
%|S_1 \cap \cdots \cap S_a| \geq n-a+1.
%\]
%This is equivalent to $\{S_1, \ldots, S_a\}$ being a circuit of
%the matroid $\H_{n,m}$, which is precisely what we wanted to show.
\end{proof}

\begin{corollary}
The matroid $\H_{n,m}$ is dual to the $(m-n+1)$-th Dilworth
truncation of the free matroid $M_m$.
\end{corollary}

\begin{proof}
This is an immediate consequence of our proof of Theorem
\ref{thm:generichyperplanes}.
\end{proof}

\section{From lines in a flag arrangement to lattice points in a simplex.}
\label{section:flags.to.simplex}

Having understood the matroid of lines in a generic hyperplane
arrangement, we proceed to study the case of complete flags. In
the following three sections, we will describe the matroid of
lines of a generic arrangement of $d$ complete flags in $\C^n$. We
start, in this section, with a short discussion of the
combinatorial setup that we will use to encode these geometric
objects. We then propose, in Section \ref{section:T_nd}, a
combinatorial definition of the matroid $\T_{n,d}$. Finally, we
will show in Section \ref{section:T_nd.is.right} that this is,
indeed, the matroid we are looking for.

Let $\E^1, \ldots, \E^d$ be $d$ generically chosen complete flags
in $\C^n$. Write
\[
\E^k = \{\{0\}=E^k_0 \subset E^k_1 \subset \cdots \subset E^k_n =
\C^n\},
\]
where $E^k_i$ is a vector space of dimension $i$.

These $d$ flags determine a line arrangement ${\bf E}_{n,d}$ in
$\C^n$ as follows. Look at all the possible intersections of the
subspaces under consideration; they are of the form
$E_{a_1,\ldots,a_d} = E^1_{a_1} \cap E^2_{a_2} \cap \cdots \cap
E^d_{a_d}$. We are interested in the one-dimensional
intersections. Since the $\E^k$s were chosen generically,
$E_{a_1,\ldots,a_d}$ has codimension $(n-a_1) + \ldots + (n-a_d)$
(or $n$ if this sum exceeds $n$). Therefore, the one-dimensional
intersections are the lines $E_{a_1,\ldots,a_d}$ for $a_1 + \cdots
+ a_d = (d-1)n+1$. There are ${n+d-2 \choose d-1}$ such lines,
corresponding to the ways of writing $n-1$ as a sum of $d$
nonnegative integers $n-a_1, \ldots, n-a_d$.

Let $T_{n,d}$ be the set of lattice points in the following
$(d-1)$-dimensional simplex in $\R^{d}$:
\[
\{\, (x_1, \ldots, x_d) \in \R^{d} \,\, | \,\, x_1 + \cdots +x_d =
n-1 \textrm{ and } x_i \geq 0 \textrm { for all } i\}.
\]
The $d$ vertices of this simplex are $(n-1,0,0,\ldots,0),
(0,n-1,0,\ldots,0), \ldots,$ $(0,0,\ldots,n-1)$.
%
%Let $\Delta_d(n)$ be the $(d-1)$-dimensional simplex in $\R^{d}$,
%whose $d$ vertices have coordinates $(n-1,0,0,\ldots,0),
%(0,n-1,0,\ldots,0), \ldots, (0,0,\ldots,n-1)$:
%\[
%\Delta_d(n) := \{\, (x_1, \ldots, x_d) \in \R^{d} \,\, | \,\, x_1
%+ \cdots +x_d = n-1 \textrm{ and } x_i \geq 0 \textrm { for all }
%i\}.
%\]
For example, $T_{n,3}$ is simply a triangular array of dots
\emph{of size $n$}; that is, with $n$ dots on each side. We will
call $T_{n,d}$ the \emph{$(d-1)$-simplex of size n}.

%We call this the simplex \emph{of size $n$}, because each edge
%contains $n$ lattice points. Let $T_{n,d}$ be the set of lattice
%points in $\Delta_d(n)$. For example, $T_{n,3}$ is simply a
%triangular array of dots \emph{of size $n$}.

It will be convenient to identify the line $E_{a_1, \ldots, a_d}$
(where $a_1+\cdots+a_d = (d-1)n+1$ and $1 \leq a_i \leq n$) with the
vector $(n-a_1, \ldots, n-a_d)$ of codimensions. This clearly gives
us a one-to-one correspondence between the set $T_{n,d}$ and the
lines in our line arrangement ${\bf E}_{n,d}$.

\begin{figure}[h]
\centering
\includegraphics[height=5.5cm]{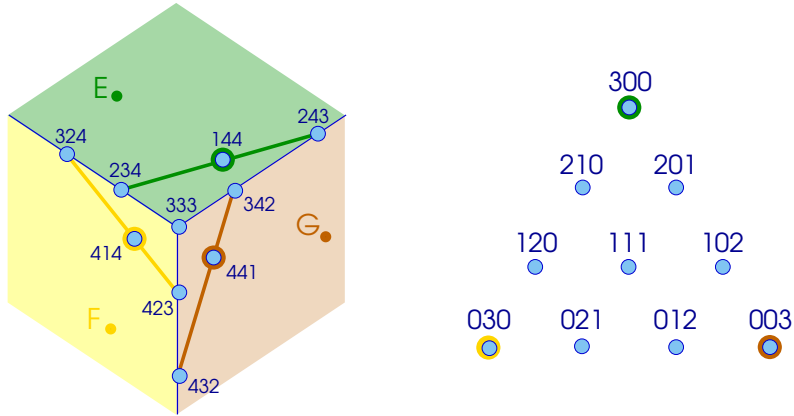}
\caption{The lines determined by three flags in $\C^4$, and the
array $T_{4,3}$.} \label{fig:flags}
\end{figure}

We illustrate this correspondence for $d=3$ and $n=4$ in Figure
\ref{fig:flags}. This picture is easier to visualize in real
projective $3$-space. Now each one of the flags $\E, \F,$ and $\G$
is represented by a point in a line in a plane. The lines in our
line arrangement are now the $10$ intersection points we see in
the picture.

We are interested in the dependence relations among the lines in
the line arrangement ${\bf E}_{n,d}$. As in the case of hyperplane
arrangements, there are several \emph{combinatorial relations}
which arise as follows. Consider a $k$-dimensional subspace
$E_{b_1, \ldots, b_d}$ with $b_1 + \cdots + b_d = (d-1)n+k$. Every
line of the form $E_{a_1, \ldots, a_d}$ with $a_i \leq b_i$ is in
this subspace, so no $k+1$ of them can be independent. The
corresponding points $(n-a_1, \ldots, n-a_d)$ are the lattice
points inside a parallel translate of $T_{k,d}$, the simplex of
size $k$, in $T_{n,d}$. In other words, in a set of independent
lines of our arrangement, we cannot have more than $k$ lines whose
corresponding dots are in a simplex of size $k$ in $T_{n,d}$.

For example, no four of the lines $E_{144}, E_{234}, E_{243},
E_{324}, E_{333},$ and $E_{342}$ are independent, because they are
in the $3$-dimensional hyperplane $E_{344}$. The dots
corresponding to these six lines form the upper $T_{3,3}$ found in
our $T_{4,3}$.

In principle, there could be other hidden dependence relations
among the lines in ${\bf E}_{n,d}$. The goal of the next two
sections is to show that this is not the case. In fact, these
combinatorial relations are the only dependence relations of the
line arrangement associated to $d$ generically chosen flags in
$\C^n$.

We will proceed as in the case of hyperplane arrangements. We will
start by showing, in Section \ref{section:T_nd}, that the
combinatorial relations do give rise to a matroid $\T_{n,d}$. In
Section \ref{section:T_nd.is.right}, we will then show that this
is, indeed, the matroid we are looking for.

%
%\begin{theorem}
%%A collection of lines in the line arrangement ${\bf E}_{n,d}$ is
%%dependent if and only if there do not exist $k+1$ of them
%The only dependence relations among the lines in the line
%arrangement ${\bf E}_{n,d}$ are the combinatorial relations.
%\end{theorem}

%Some partial evidence for this claim is provided by Theorem
%\ref{thm:matroid}, which we now present.

\section{A matroid on the lattice points in a regular \\simplex.}
\label{section:T_nd}

%The lines in the line arrangement determined by $d$ generic flags
%in $\C^n$ correspond to the dots in $T_{n,d}$. We wish to show
%that a subset of the line arrangement is independent if and only
%if the corresponding subset of $T_{n,d}$ is in the collection
%$\I_{n,d}$.
%
%Before we prove this, we will show that $\I_{n,d}$ is, indeed, the
%collection of independent sets of a matroid. This will provide a
%better understanding of the combinatorial structure of the
%dependence relations we are dealing with.

In this section, we show that the combinatorial dependence
relations defined in Section \ref{section:flags.to.simplex} do
determine a matroid.

\begin{theorem}\label{thm:matroid}
Let $\I_{n,d}$ be the collection of subsets $I$ of $T_{n,d}$ such
that every parallel translate of $T_{k,d}$ contains at most $k$
points of $I$, for every $k \leq n$.

Then $\I_{n,d}$ is the collection of independent sets of a matroid
$\T_{n,d}$ on the ground set $T_{n,d}$.
\end{theorem}

\begin{figure}[h]
\centering
\includegraphics[height=3cm]{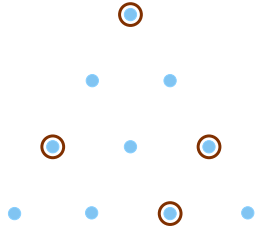}
\caption{The array $T_{4,3}$ and a set in $\I_{4,3}$.}
\label{fig:T34}
\end{figure}

We will call a parallel translate of $T_{k,d}$ a \emph{simplex of
size $k$}. As an example, $T_{n,3}$ is a triangular array of dots
\emph{of size $n$}. The collection $\I_{n,3}$ consists of those
subsets $I$ of the array $T_{n,3}$ such that no triangle of size
$m$ contains more than $m$ points of $I$. Figure \ref{fig:T34}
shows the array $T_{4,3}$, and a set in $\I_{4,3}$.

\begin{proof}[Proof of Theorem \ref{thm:matroid}.]
We need to verify the three axioms for the collection of
independent sets of a matroid:

\noindent{\bf (I1)} The empty set is in $\I_{n,d}$.

\noindent{\bf (I2)} If $I$ is in $\I_{n,d}$ and $I' \subseteq I$,
then $I'$ is also in $\I_{n,d}$.

\noindent{\bf (I3)} If $I$ and $J$ are in $\I_{n,d}$ and $|I| <
|J|$, then there is an element $e$ in $J-I$ such that $I \cup e$
is in $\I_{n,d}$.

The first two axioms are satisfied trivially; let us focus on the
third one. Proceed by contradiction. Let $J-I = \{e_1, \ldots,
e_k\}$. We know that every simplex of size $a$ contains at most
$a$ points of $I$. When we try to add $e_h$ to $I$ while
preserving this condition, only one thing can stop us: a simplex
$T_h$ of size $a_h$ which already contains $a_h$ points of $I$,
and also contains $e_h$.

Say that a simplex of size $a$ is \emph{$I$-saturated} if it
contains exactly $a$ points of $I$. We have found $I$-saturated
simplices $T_1, \ldots, T_k$ which contain $e_1, \ldots, e_k$,
respectively.

Now we use the following lemma, which we will prove in a moment.

\begin{lemma}\label{lemma:saturated}
Let $T$ and $T'$ be two $I$-saturated simplices with $T \cap T'
\neq \emptyset$. Let $T \vee T'$ be the smallest simplex
containing $T$ and $T'$. Then the simplices $T \cap T'$ and $T
\vee T'$ are also $I$-saturated.
\end{lemma}

If two of our $I$-saturated simplices $T_g$ and $T_h$ are
different and have a non-empty intersection, we can replace them
both by $T_g \vee T_h$. By Lemma \ref{lemma:saturated}, this is a
larger $I$-saturated simplex, and it still contains $e_g$ and
$e_h$. We can continue in this way, until we obtain $I$-saturated
simplices $T_1, \ldots, T_k$ containing $e_1, \ldots, e_k$
\emph{which are pairwise disjoint} (though possibly repeated).

Let $U_1, \ldots, U_l$ be this collection of $I$-saturated
simplices, now listed without repetitions. Let $U_r$ have size
$s_r$, and say it contains $i_r$ elements of $I-J$, $j_r$ elements
of $J-I$, and $k_r$ elements of $I \cap J$.

We know that $U_r$ is $I$-saturated, so $s_r = i_r + k_r.$ We also
know that $J$ is in $\I_{n,d}$, so $s_r \geq j_r + k_r$.
Therefore, $i_r \geq j_r$ for each $r$.

Now, the $U_r$s are pairwise disjoint, so $\sum i_r \leq |I-J|$
and $\sum j_r \leq |J-I|$. But in fact, we know that every element
of $J-I$ is in some $U_r$, so we actually have the equality $\sum
j_r = |J-I|$. Therefore we have
\[
|J-I| = \sum j_r \leq \sum i_r \leq |I-J|.
\]
This contradicts our assumption that $|I| < |J|$, and Theorem
\ref{thm:matroid} follows.
\end{proof}

\emph{Proof of Lemma \ref{lemma:saturated}.} First we show that
$\s(T \cap T') + \s(T \vee T') = \s(T) + \s(T')$. Each simplex is
a parallel translate of some $T_{k,d}$; its vertices are given by
$(a_1+k-1,a_2,\ldots,a_d), \ldots, (a_1,a_2, \ldots, a_d+k-1)$ for
some $a_1,\ldots,a_d$ such that $\sum a_i = n-k$. We denote this
simplex by $T_{a_1,\ldots,a_d}$; its size is $k=n-\sum a_i$. It
consists of the points $(x_1,\ldots,x_d)$ with $x_i \geq a_i$ for
each $i$, and $\sum x_i = n-1$. Therefore, $T_{a_1,\ldots,a_d}
\subseteq T_{A_1,\ldots,A_d}$ if and only if $a_i \geq A_i$ for
each $i$.

It follows that if we let $T=T_{a_1,\ldots,a_d}$ and
$T'=T_{a_1',\ldots,a_d'}$, then we have:
\begin{eqnarray*}
T \cap T' &=& T_{\max(a_1,a_1'), \ldots, \max(a_d,a_d')} \\
T \vee T' &=& T_{\min(a_1,a_1'),\ldots,\min(a_d,a_d')}.
\end{eqnarray*}
So $\s(T \cap T') + \s(T \vee T') = (n-\sum \max(a_i,a_i')) + (n -
\sum \min(a_i,a_i'))$ and $\s(T)+\s(T')= (n-\sum a_i )+(n-\sum
a_i')$. These are equal since $\max(a,a')+\min(a,a')=a+a'$ for any
$a,a' \in \R$.

%Now, it is clear that $|I \cap (T \cap T')| + |I \cap (T \vee T')|
%\geq |I \cap T| + |I \cap T| = \s(T) +\s(T') = \s(T \cap T') +
%\s(T \vee T')$. Since $I$ is in $\I_{n,d}$, this implies $|I \cap
%(T \cap T')| = \s(T \cap T')$ and $|I \cap (T \vee T')| = \s(T
%\vee T')$, as desired. $\qed$

We know that $T$ and $T'$ contain $\s(T)$ and $\s(T')$ points of
$I$, respectively. If $T \cap T'$ and $T \vee T'$ contain $x$ and
$y$ points of $I$, we have that $x+y \geq \s(T)+\s(T')$, so $x+y
\geq \s(T \cap T') + \s(T \vee T')$. But $I$ is in $\I_{n,d}$, so
$x \leq \s(T \cap T')$ and $y \leq \s(T \vee T')$. This can only
happen if equality holds, and $T \cap T'$ and $T \vee T'$ are
$I$-saturated. $\qed$
%
%\medskip
%
%We can say more about the structure of the matroids $\H_{n,m}$ and
%$\T_{n,d}$. Given a set $E$ and a function $f:2^E \rightarrow \N$
%which is increasing ($f(X) \leq f(Y)$ for $X \subseteq Y$) and
%submodular ($f(X)+f(Y) \geq f(X \cup Y) + f(X \cap Y)$), there
%exists a matroid $M(f)$ on $E$ whose circuits are the minimal sets
%$S \subseteq E$ such that $|S|
%> f(S)$. This matroid is denoted $M(f)$, and called the
%\emph{polymatroid} of the function $f$. \cite[Section 12.1]{Oxley}
%%
%%For example, consider the function $f$ which assigns, to a
%%collection $S$ of $(n-1)$-subsets of $[m]$, the number $f(S) := n
%%- \bigl|\bigcap_{A \in S} A\bigr|$. This function is easily seen
%%to be increasing and submodular, and it is clear from the
%%definition that the matroid $\H_{n,m}$ of Proposition
%%\ref{prop:H_mn} is the polymatroid $M(f)$.
%
%\begin{corollary}
%\begin{enumerate}
%\item The matroid $\H_{n,m}$ is the polymatroid $M(f)$, where $f$
%is the increasing and submodular function that assigns the number
%\[
%f(S) := n - \bigl|\hspace{-.1cm}\bigcap_{A \in S} A\bigr|
%\]
%to a collection $S$ of $(n-1)$-subsets of $[m]$. \item  The
%matroid $\T_{n,d}$ is the polymatroid $M(g)$, where $g$ is the
%increasing and submodular function which assigns the number
%\[
%g(S) = n - \min_{a \in S} a_1 - \cdots - \min_{a \in S} a_d
%\]
%to a collection of points $S$ in $T_{n,d}$.
%\end{enumerate}
%\end{corollary}
%
%\begin{proof}
%These statements follow directly from Proposition \ref{prop:H_mn}
%and the proof of Theorem \ref{thm:matroid}.
%\end{proof}

\section{This is the right matroid.} \label{section:T_nd.is.right}

We now show that the matroid $\T_{n,d}$ of Section
\ref{section:T_nd} is, indeed, the matroid that arises from
intersecting $d$ flags in $\C^n$ which are generic enough.
%
%As before, we consider $d$ generically chosen complete flags in
%$\C^n$, and the ${n+d-2 \choose d-1}$ lines they determine, which
%are of the form $E_{a_1, \ldots, a_d} = E^1_{a_1} \cap \cdots \cap
%E^d_{a_d}$ for $1 \leq a_i \leq n$ and $a_1 + \cdots + a_d =
%(d-1)n+1$.

\begin{theorem}\label{thm:genericflags}
If $d$ complete flags $\E^1, \ldots, \E^d$ in $\C^n$ are generic
enough, then the matroid of the ${n+d-2 \choose d-1}$ lines
$E_{a_1,\ldots,a_d}$ is isomorphic to $\T_{n,d}$.
\end{theorem}

\begin{proof}
As mentioned in Section \ref{section:flags.to.simplex}, the
one-dimensional intersections of the $\E^i$s satisfy the following
combinatorial relations: each $k$ dimensional subspace
$E_{b_1\ldots b_d}$ with $b_1 + \cdots +b_d = (d-1)n+k$, contains
the lines $E_{a_1\ldots a_d}$ with $a_i \leq b_i$; therefore, it
is impossible for $k+1$ of these lines to be independent. The
subspace $E_{b_1\ldots b_d}$ corresponds to the simplex of dots
which is labelled $T_{n-b_1, \ldots, n-b_d}$, and has size $n -
\sum(n-b_i) = k$. The lines $E_{a_1\ldots a_d}$ with $a_i \leq
b_i$ correspond precisely the dots in this copy of $T_{k,d}$. So
these ``combinatorial relations" are precisely the dependence
relations of $\T_{n,d}$.

Now we need to show that, if the flags are ``generic enough",
these are the only linear relations among these lines. It is
enough to construct one set of flags which satisfies no other
relations.

Consider a set $\H$ of $d(n-1)$ hyperplanes $H^i_j$ in $\C^n$ (for
$1 \leq i \leq d$ and $1 \leq j \leq n-1$) which are generic in
the sense of Theorem \ref{thm:generichyperplanes}, so the only
dependence relations among their one-dimensional intersections are
the combinatorial ones. Now, for $i = 1, \ldots, d$, define the
flag $\E^i$ by:
\begin{eqnarray*}
E^i_{n-1} & = & H^i_{n-1} \\
E^i_{n-2} & = & H^i_{n-1} \cap H^i_{n-2} \\
& \vdots & \\
E^i_1 & = & H^i_{n-1} \cap H^i_{n-2} \cap \cdots \cap H^i_1,
\end{eqnarray*}
We will show that these $d$ flags are generic enough; in other
words, the matroid of their one-dimensional intersections is
$\T_{n,d}$.

Let us assume that a set $S$ of one-dimensional intersections of the
$\E^i$s is dependent. Since each line in $S$ is a one-dimensional
intersection of the hyperplanes $H^i_j$, we can apply Theorem
\ref{thm:generichyperplanes}. It tells us that for some $t$ we can
find $t+1$ lines in $S$
%(say
%$E_{a_1, \ldots, a_d}, \ldots, E_{z_1, \ldots, z_d}$)
and a set $T$ of $n-t$ hyperplanes $H^i_j$ which contain all of
them.

Our $t+1$ lines are of the form
\begin{eqnarray*}
E_{a_1, \ldots, a_d} &=& E^1_{a_1} \cap \cdots \cap E^d_{a_d} \\
&=& (H^1_{n-1} \cap \cdots \cap H^1_{a_1}) \cap \cdots \cap
(H^d_{n-1} \cap \cdots \cap H^d_{a_d}).
\end{eqnarray*}
Therefore, if a hyperplane $H^i_j$ contains them, so does $H^i_k$
for any $k > j$. Let us add all such hyperplanes to our set $T$,
to obtain the set
\[
U = \{H^1_{n-1}, \ldots, H^1_{b_1}, \ldots, H^d_{n-1}, \ldots,
H^d_{b_d}\},
\]
where $b_i$ is the smallest $j$ for which $H^i_j$ is in $T$. The
set $U$ contains $\sum (n-b_i)$ hyperplanes, so $ \sum (n-b_i)
\geq n-t.$

Each one of our $t+1$ lines is contained in each of the
hyperplanes in $U$, and therefore in their intersection
\[
\bigcap_{H^i_j \in U} H^i_j = E_{b_1,\ldots,b_d},
\]
which has dimension $n - \sum(n-b_i) \leq t$.

So, actually, the dependence of the set $S$ is a consequence of
one of the combinatorial dependence relations present in
$\T_{n,d}$. The desired result follows.
\end{proof}
%
%\begin{figure}[h]
%\centering
%\includegraphics[height=5cm]{flags3withbasis}
%\caption{A basis of $T_{3,4}$ \textsc{Fix labels.}}
%\label{fig:flags}
%\end{figure}

With Theorem \ref{thm:genericflags} in mind, we will say that the
complete flags $\E^1, \ldots, \E^d$ in $\C^n$ are
\emph{matroid-generic} if the matroid of the ${n+d-2 \choose d-1}$
lines $E_{a_1,\ldots,a_d}$ is isomorphic to $\T_{n,d}$.

We conclude this section by showing that the one-dimensional
intersections of matroid-generic flag arrangements are completely
characterized by their combinatorial properties.

\begin{proposition} \label{prop:matroid.to.arrangement.}
If a line arrangement $\L$ in $\C^n$ has matroid $\T_{n,d}$, then
it can be realized as the arrangement of one-dimensional
intersections of $d$ complete flags in $\C^n$.
\end{proposition}

\begin{proof}
To make the notation clearer, let us give the proof for $d=3$, which
generalizes trivially to larger values of $d$. Denote the lines in
$\L$ by $L_{rst}$ for $r+s+t=2n+1$. Consider the three flags $\E,
\F$ and $\G$ given by
\begin{eqnarray*}
E_i &=& \textrm{span} \{L_{rst} \,\, | \,\, r \leq i\} \\
F_i &=& \textrm{span} \{L_{rst} \,\, | \,\, s \leq i\} \\
G_i &=& \textrm{span} \{L_{rst} \,\, | \,\, t \leq i\}
\end{eqnarray*}
for $0 \leq i \leq n$. Compare this with Figure \ref{fig:flags} in
the case $n=4$. The subspace $E_i$, for example, is the span of the
lines corresponding to the first $i$ rows of the triangle.

Since $\L$ is a representation of the matroid $\T_{n,3}$, the
dimensions of $E_i$, $F_i$, and $G_i$ are equal to $i$, which is the
rank of the corresponding sets (copies of $T_{i,3}$) in $\T_{n,3}$.

We now claim that the line arrangement corresponding to $\E, \F$ and
$\G$ is precisely $\L$. This amounts to showing that $E_i \cap F_j
\cap G_k = L_{ijk}$ for $i+j+k=2n+1$. We know that $L_{ijk}$ is in
$E_i$, $F_j$, and $G_k$ by definition, so we simply need to show
that $\dim(E_i \cap F_j \cap G_k) = 1$.

Assume $\dim(E_i \cap F_j \cap G_k) \geq 2$. Consider the sequence
of subspaces:
%\begin{eqnarray*}
%E_i \cap F_j \cap G_k & \subseteq & E_{i+1} \cap F_j \cap G_k & \subseteq & \cdots & \subseteq & E_n \cap F_j \cap G_k \\%
%& \subseteq & E_n \cap F_{j+1} \cap G_k & \subseteq & \cdots & \subseteq & E_n \cap F_n \cap G_k \\
%& \subseteq & E_n \cap F_n \cap G_{k+1} & \subseteq & \cdots &
%\subseteq & E_n \cap F_n \cap G_n
%\end{eqnarray*}

\begin{tabular}{cccccccccc}
$E_i \cap F_j \cap G_k$ & $\subseteq$ & $E_{i+1} \cap F_j \cap G_k$
& $\subseteq$ & $\cdots$ & $\subseteq$ & $E_n \cap F_j \cap G_k$ &
$\subseteq$ & \\
& $\subseteq$ & $E_n \cap F_{j+1} \cap G_k$ & $\subseteq$ & $\cdots$
& $\subseteq$ & $E_n \cap F_n \cap G_k$ &
$\subseteq$ & \\
& $\subseteq$ & $E_n \cap F_n \cap G_{k+1}$ & $\subseteq$ & $\cdots$
& $\subseteq$ & $E_n \cap F_n \cap G_n$. &
\end{tabular}

There are $1+(n-i)+(n-j)+(n-k) = n$ subspaces on this list; the
first one has dimension at least $2$, and the last one has dimension
$n$. By the pigeonhole principle, two consecutive subspaces on this
list must have the same dimension. Since one is contained in the
other, these two subspaces must actually be equal. So assume that
$E_{a-1} \cap F_b \cap G_c = E_a \cap F_b \cap G_c$; a similar
argument will work in the other cases.

Now, we have $a+b+c > i+j+k = 2n+1$, so we can find positive
integers $\beta \leq b$ and $\gamma \leq c$ such that
$a+\beta+\gamma = 2n+1$. Then $L_{a\beta\gamma}$ is a line which, by
definition, is in $E_a$, $F_b$ and $G_c$. It follows that
\[
L_{a\beta\gamma} \in E_a \cap F_b \cap G_c = E_{a-1} \cap F_b \cap
G_c \subseteq E_{a-1}.
\]
This implies that $L_{a\beta\gamma}$ is dependent on $\{L_{rst} \,\,
| \,\, r \leq a-1\}$, which is impossible since $\L$ represents the
matroid $\T_{n,3}$. We have reached a contradiction, which implies
that $\dim(E_i \cap F_j \cap G_k) = 1$ and  therefore $E_i \cap F_j
\cap G_k = L_{ijk}$.

%Assume  $\dim(E_i \cap F_j \cap G_k) \geq 2$. Since $\L$ has rank
%$\rk(\T_{n,3})=n$, it spans $\C^n$; therefore we can find some
%other $L_{rst}$ with $r+s+t=2n+1$ which is also in $E_i \cap F_j
%\cap G_k$. Since $r+s+t=i+j+k$, we can assume, without loss of
%generality, that $r>i$. But then, clearly, the set $(r,s,t) \cup
%\{(a,b,c) | a \leq i\}$ has rank $i+1$ in $\T_{n,3}$. Since $\L$
%represents $\T_{n,3}$, this implies that $L_{rst}$ is not in
%$E_i$, a contradiction.

It follows that $\L$ is the line arrangement determined by flags
$\E, \F$ and $\G$, as we wished to show.
\end{proof}

\section{Rhombus tilings of holey triangles and the matroid
$\T_{n,3}$.}\label{section:rhombus.tilings}

Let us change the subject for a moment.

\begin{figure}[h]
\centering
\includegraphics[height=3cm]{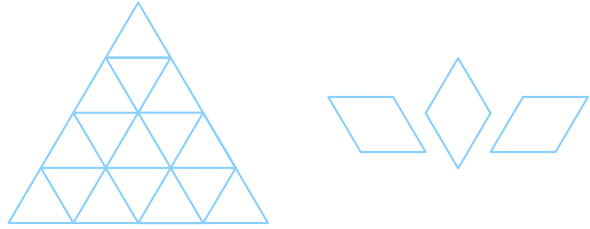}
\caption{$T(4)$ and the three rhombus tiles.} \label{fig:T(4)}
\end{figure}

Let $T(n)$ be an equilateral triangle with side length $n$.
Suppose we wanted to tile $T(n)$ using unit rhombi with angles
equal to $60^{\circ}$ and $120^{\circ}$. It is easy to see that
this task is impossible, for the following reason. Cut $T(n)$ into
$n^2$ unit equilateral triangles, as illustrated in Figure
\ref{fig:T(4)}; $n(n+1)/2$ of these triangles point upward, and
$n(n-1)/2$ of them point downward. Since a rhombus always covers
one upward and one downward triangle, we cannot use them to tile
$T(n)$.

\begin{figure}[h]
\centering
\includegraphics[height=3cm]{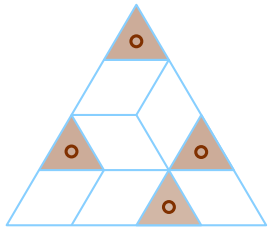}
\caption{A tiling of a holey $T(4)$.} \label{fig:tiling}
\end{figure}

Suppose then that we make $n$ holes in the triangle $T(n)$ by
cutting out $n$ of the upward triangles. Now we have an equal
number of upward and downward triangles, and it may or may not be
possible to tile the remaining shape with rhombi. Figure
\ref{fig:tiling} shows a tiling of one such \emph{holey triangle}.

The main question we address in this section is the following:
\begin{question}\label{question}
Given $n$ holes in $T(n)$, is there a simple criterion to
determine whether there exists a rhombus tiling of the holey
triangle that remains?
\end{question}

A rhombus tiling is equivalent to a perfect matching between the
upward triangles and the downward triangles. Hall's theorem then
gives us an answer to Question \ref{question}: It is necessary and
sufficient that any $k$ downward triangles have a total of at
least $k$ upward triangles to match to.

%We can represent upward and downward triangles by black and white
%vertices, respectively, and join a black and a white vertex when
%the corresponding triangles are adjacent. Finding a rhombus tiling
%of the holey triangle is equivalent to constructing a matching of
%the dual bipartite graph, and Hall's theorem gives us an answer to
%Question \ref{question}: any $k$ downward triangles must have a
%total of at least $k$ upward triangles to match to.

However, the geometry of $T(n)$ allows us to give a simpler
criterion. Furthermore, this criterion reveals an unexpected
connection between these rhombus tilings and the line arrangement
determined by $3$ generically chosen flags in $\C^n$. Notice that
the upward triangles in $T(n)$ can be identified with the dots of
$T_{n,3}$.

\begin{theorem}\label{thm:tilings}
Let $S$ be a set of $n$ holes in $T(n)$. The triangle $T(n)$ with
holes at $S$ can be tiled with rhombi if and only if the locations
of the holes constitute a basis for the matroid $\T_{n,3}$;
\emph{i.e.}, if and only if every $T(k)$ in $T(n)$ contains at
most $k$ holes of $S$, for all $k \leq n$.
\end{theorem}

\begin{proof}
First suppose that we have a tiling of the holey triangle, and
consider any triangle $T(k)$ in $T(n)$. Consider all the tiles
which contain one or two triangles of that $T(k)$, and let $R$ be
the holey region that these tiles cover. Since the boundary of
$T(k)$ consists of upward triangles, the region $R$ is just $T(k)$
with some downward triangles glued to its boundary.

If $T(k)$ had more than $k$ holes, it would have fewer than
$k(k-1)/2$ upward triangles, and so would $R$. However, $R$ has at
least the $k(k-1)/2$ downward triangles of $T(k)$. That makes it
impossible to tile the region $R$, which contradicts its
definition. This proves the forward direction.

\medskip

Now let $S$ be a set of $n$ holes in $T(n)$ such that every $T(k)$
contains at most $k$ holes. Equivalently, think of $S$ as a basis
of the matroid $\T_{n,3}$. We construct a tiling of the resulting
holey triangle by induction on $n$. The case $n=1$ is trivial, so
assume $n \geq 2$.

Within that induction, we induct on the number of holes of $S$ in
the bottom row of $T(n)$. Since the $T(n-1)$ of the top $n-1$ rows
contains at most $n-1$ holes, there is at least one hole in the
bottom row.

If there is exactly one hole in the bottom row, then the tiling of
the bottom row is forced upon us, and the top $T(n-1)$ can be
tiled by induction. Now assume that there are at least two holes
in the bottom row; call the two leftmost holes $x$ and $y$ in that
order. Consider the upward triangles in the second to last row
which are between $x$ and $y$; label them $a_1, \ldots, a_t$. This
is illustrated in an example in the top left panel of Figure
\ref{fig:slidenice}. Here $a_1, a_2, a_3$ and $a_4$ are shaded
lightly, and $a_1$ is also a hole.

\begin{figure}[h]
\centering
\includegraphics[height=4.5cm]{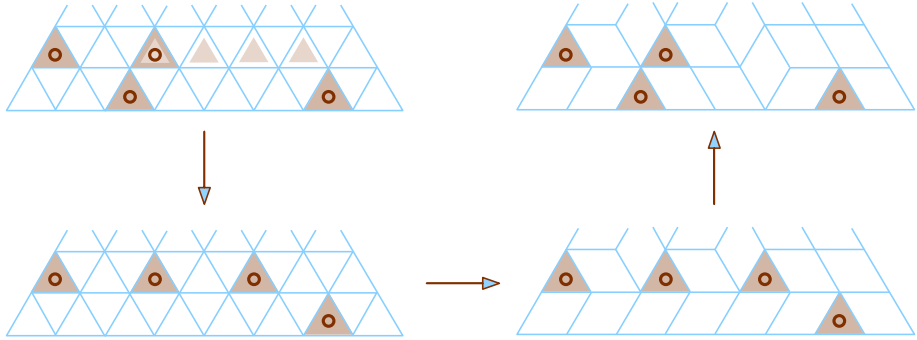}
\caption{Sliding the hole from $a_i$ to $x$.} \label{fig:slidenice}
\end{figure}

We claim that we can exchange the hole $x$ for one of the holes
$a_i$, so that the set of holes $S-x\cup a_i$ is also a basis of
$\T_{n,3}$. Notice that this $a_i$ cannot be in $S$. Assume that no
such $a_i$ exists. Then each $a_i$ must be in a triangle $T_i$ which
is $(S-x)$-saturated.\footnote{As in Section \ref{section:T_nd}, if
$A$ is a set of holes, we say that an upward triangle of size $k$ is
\emph{$A$-saturated} if it contains $k$ holes of $A$.} If $a_i$ is
in $S$, then $T_i=a_i$. The triangle $y$ is also trivially
$(S-x)$-saturated. Observe that Lemma \ref{lemma:saturated} also
applies in this setting, even when two saturated triangles intersect
only on their boundary. We can use this successively to obtain an
$(S-x)$-saturated triangle containing $a_1, \ldots, a_t$, and $y$.
But that triangle will also contain $x$, so it will contain more
holes of $S$ than it is allowed.

So let $a_i$ be such that $S-x\cup a_i$ is a basis of $\T_{n,3}$.
For instance, in the first step of Figure \ref{fig:slidenice}, $x$
is exchanged for $a_3$. Notice that $S-x\cup a_i$ contains fewer
holes in the bottom row than $S$ does. By the induction
hypothesis, we can tile the $T(n)$ with holes at $S-x \cup a_i$,
as shown in the second step of Figure \ref{fig:slidenice}. The
bottom row of this tiling is frozen from left to right until it
reaches $y$. Therefore, we can slide the hole from $a_i$ back to
$x$ in the obvious way, by reversing the tiles in the bottom row
between $x$ and $a_i$. This is illustrated in the last step of
Figure \ref{fig:slidenice}. We are left with a tiling with holes
at $S$, as desired.
\end{proof}

%
%\begin{corollary}\label{thm:tilings}
%The possible locations of $n$ holes for which a rhombus tiling of
%the holey triangle $T(n)$ exists correspond to the bases of the
%matroid $\T_{n,3}$.
%\end{corollary}
%
%\begin{proof}
%This is just a restatement of Theorem \ref{thm:tilings}.
%\end{proof}

Theorem \ref{thm:tilings} allows us to say more about the
structure of the matroid $\T_{n,3}$. We first remind the reader of
the definition of two important families of matroids, called
\emph{transversal} and \emph{cotransversal} matroids. For more
information, we refer the reader to \cite{Ardila, Oxley}.

%Let $M$ be a matroid on a finite ground set $S$, whose collection
%of bases is $\B$. Define $\B^* = \{S-B \,\, | \,\, B \in \B\}$ be
%the collection of complements of the bases of $B$. Then $\B^*$ is
%also the collection of bases of a matroid on the ground set $S$,
%which is called the \emph{dual matroid} of $M$, and denoted $M^*$.

Let $S$ be a finite set, and let $A_1, \ldots, A_r$ be subsets of
$S$. A \emph{transversal} of $(A_1, \ldots, A_r)$, also known as a
\emph{system of distinct representatives}, is a subset $\{e_1,
\ldots, e_r\}$ of $S$ such that $e_i$ is in $A_i$ for each $i$,
and the $e_i$s are distinct. The transversals of $(A_1, \ldots,
A_r)$ are the bases of a matroid on $S$. Such a matroid is called
a \emph{transversal matroid}, and $(A_1, \ldots, A_r)$ is called a
\emph{presentation} of the matroid.

Let $G$ be a directed graph with vertex set $V$, and let $A=\{v_1,
\ldots, v_r\}$ be a subset of $V$. We say that an $r$-subset $B$
of $V$ \emph{can be linked to $A$} if there exist $r$
vertex-disjoint directed paths whose initial vertex is in $B$ and
whose final vertex is in $A$. We will call these $r$ paths a
\emph{routing} from $B$ to $A$. The collection of $r$-subsets
which can be linked to $A$ are the bases of a matroid denoted
$L(G,A)$. Such a matroid is called a \emph{cotransversal} matroid
or a \emph{strict gammoid}. It is a nontrivial fact that these
matroids are precisely the duals of the transversal matroids
\cite{Ardila, Oxley}, and we will call them \emph{cotransversal}.

\begin{theorem}\label{thm:transversal}
The matroid $\T_{n,3}$ is cotransversal.
\end{theorem}

\begin{proof}[First proof.]
We prove that $\T_{n,3}^*$ is transversal. We can think of the
ground set of $\T_{n,3}$ as the set of upward triangles in $T(n)$.
By Theorem \ref{thm:tilings}, a basis of $\T_{n,3}$ is a set of
$n$ holes for which the resulting holey triangle can be tiled; its
complement is the set of ${n \choose 2}$ upward triangles which
share a tile with one of the ${n \choose 2}$ downward triangles.

Number the downward triangles $1, 2, \ldots, {n \choose 2}$ from
top to bottom and left to right. Then a  tiling of the complement of a basis of
$\T_{n,3}$ is nothing but a transversal of $(A_1, \ldots, A_{n
\choose 2})$, where $A_i$ is the set of three upward triangles
which are adjacent to downward triangle $i$. This completes the
proof.
\end{proof}

\begin{figure}[h]
\centering
\includegraphics[height=4cm]{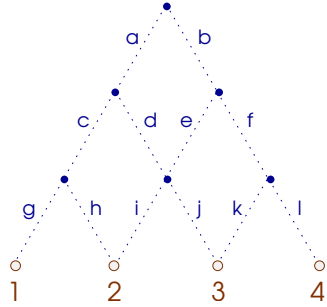}
\caption{The graph $G_4$.} \label{fig:g4}
\end{figure}

\begin{proof}[Second proof.]
We prove that $\T_{n,3}$ is cotransversal. Let $G_n$ be the
directed graph whose set of vertices is the triangular array
$T_{n,3}$, where each dot not on the bottom row is connected to
the two dots directly below it. Label the dots on the bottom row
$1,2,\ldots, n$. Figure \ref{fig:g4} shows $G_4$; all the edges of
the graph point down.

There is a well-known trick which allows us to view rhombus tilings
of the holey triangle $T(n)$ as routings in $G_n$. It works as
follows: The midpoints of the possible horizontal edges of a tiling
form a copy of the graph $G_{n}$. Given a tiling of a holey $T(n)$,
join two vertices of $G_{n}$ if they are on opposite edges of the
same tile; this gives the desired routing of $G_{n}$. This
correspondence is best understood in an example; see Figure
\ref{fig:tilingwithpaths}.

\begin{figure}[h]
\centering
\includegraphics[height=4.5cm]{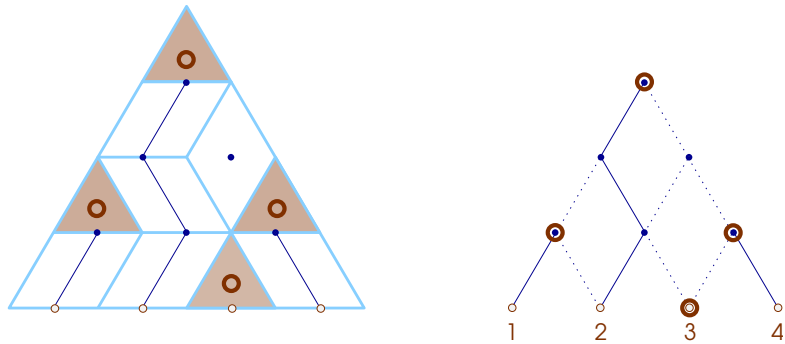}
\caption{A tiling of a holey $T(4)$ and the corresponding routing of
$G_4$.} \label{fig:tilingwithpaths}
\end{figure}

It is easy to check that this is a bijection between the rhombus
tilings of the holey triangles of size $n$, and the routings in the
graph $G_n$ which start anywhere and end at vertices $1,2,\ldots,
n$. Also, the holes of the holey triangle correspond to the starting
points of the $n$ paths in the routing. From Theorem
\ref{thm:tilings}, it follows that $\T_{n,3}$ is the cotransversal
matroid $L(G_n, [n])$.
\end{proof}

\begin{theorem}\label{thm:represent.cotransversal}
Assign algebraically independent weights to the edges of
$G_n$.\footnote{Integer weights which increase extremely quickly
will also work.} For each dot $D$ in the triangular array
$T_{n,3}$ and each $1 \leq i \leq n$, let $v_{D,i}$ be the sum of
the weights of all paths\footnote{The weight of a path is defined
to be the product of the weights of its edges.} from dot $D$ to
dot $i$ on the bottom row.

Then the \emph{path vectors} $v_D = (v_{D,1}, \ldots, v_{D,n})$
are a geometric representation of the matroid $T_{n,3}$.
\end{theorem}

For example, the top dot of $T_{4,3}$ in Figure \ref{fig:g4} would
be assigned the \emph{path vector} $(acg, ach+adi+bei,
adj+bej+bfk, bfl)$. Similarly, focusing our attention on the top
three rows, the representation we obtain for the matroid
$\T_{3,3}$ is given by the columns of the following matrix:
\[
\left(
\begin{array}{cccccc}
1 & 0 & 0 & c & 0 & ac \\
0 & 1 & 0 & d & e & ad+be \\
0 & 0 & 1 & 0 & f & bf
\end{array}
\right)
\]

\begin{proof}[Proof of Theorem \ref{thm:represent.cotransversal}.]
By the Lindstr\"om-Gessel-Viennot lemma \cite{Gessel, Karlin,
Lindstrom, Mason2}, the determinant of the matrix with columns
$v_{D_1}, \ldots, v_{D_n}$ is equal to the signed sum of the
routings from $\{D_1, \ldots, D_n\}$ to $\{1,\ldots,n\}$. The sign
of a routing is the sign of the permutation of $S_n$ which matches
the starting points and the ending points of the $n$ paths. Since
we chose the weights to be algebraically independent, this sum is
zero if and only if it is empty.

Therefore, $v_{D_1}, \ldots, v_{D_n}$ are independent if and only
if there exists a routing from $\{D_1, \ldots, D_n\}$ to
$\{1,\ldots,n\}$; that is, if and only if $\{D_1, \ldots, D_n\}$
are a basis of $L(G_n, [n])$.
\end{proof}

It is worth pointing out that Lindstr\"om's original motivation for
the discovery of the Lindstr\"om-Gessel-Viennot lemma was to explain
Mason's construction of a geometric representation of an arbitrary
cotransversal matroid \cite{Lindstrom, Mason}. Theorem
\ref{thm:represent.cotransversal} and its proof are special cases of
their more general argument; we have included them for completeness.

The very simple and explicit representation of $\T_{n,3}$ of
Theorem \ref{thm:represent.cotransversal} will be shown in Section
\ref{section:schubert} to have an unexpected consequence in the
Schubert calculus: it provides us with a reasonably efficient
method for computing Schubert structure constants in the flag
manifold.

\section{Fine mixed subdivisions of $n\Delta_{d-1}$
and triangulations of $\Delta_{n-1} \times
\Delta_{d-1}$.}\label{section:fine.mixed.subdivisions}

The surprising relationship between the geometry of three flags in
$\C^n$ and the rhombus tilings of holey triangles is useful to us
in two ways: it explains the structure of the matroid $\T_{n,3}$,
and it clarifies the conditions for a rhombus tiling of such a
region to exist. We now investigate a similar connection between
the geometry of $d$ flags in $\C^n$, and certain
$(d-1)$-dimensional analogs of these tilings, known as \emph{fine
mixed subdivisions} of $n \Delta_{d-1}$.

The fine mixed subdivisions of $n \Delta_{d-1}$ are in one-to-one
correspondence with the triangulations of the polytope
$\Delta_{n-1} \times \Delta_{d-1}$. The triangulations of a
product of two simplices are fundamental objects, which have been
studied from many different points of view. They are of
independent interest \cite{Babson, Bayer, Gelfand}, and have been
used as a building block for finding efficient triangulations of
high dimensional cubes \cite{Haiman, Orden} and disconnected
flip-graphs \cite{Santosflips, Santostoric}. They also arise very
naturally in connection with tropical geometry \cite{Develin},
transportation problems, and Segre embeddings \cite{Sturmfels}. In
the following two sections, we provide evidence that
triangulations of $\Delta_{n-1} \times \Delta_{d-1}$ are also
closely connected to the geometry of $d$ flags in $\C^n$, and that
their study can be regarded as a study of \emph{tropical oriented
matroids}.

% plays a crucial role
%in our proof of Theorem \ref{thm:3generic}. In our attempt to
%prove the analogous Conjecture \ref{conj:dgeneric} on $d$ flags in
%$\C^n$, it is natural to look for a $(d-1)$-dimensional
%generalization of our results on rhombus tilings. In this section,
%we present some promising partial results in that direction.

\medskip

Instead of thinking of rhombus tilings of a holey triangle, it
will be slightly more convenient to think of them as \emph{lozenge
tilings} of the triangle: these are the tilings of the triangle
using unit rhombi and upward unit triangles. A good
high-dimensional analogue of the lozenge tilings of the triangle
$n \Delta_2$ are the \emph{fine mixed subdivisions} of the simplex
$n \Delta_{d-1}$; we briefly recall their definition.

The \emph{Minkowski sum} of polytopes $P_1, \ldots, P_k$ in
$\R^m$, is:
\[
P = P_1 + \cdots + P_k := \{p_1+\cdots+p_k \, | \, p_1 \in P_1,
\ldots, p_k \in P_k \}.
\]
We are interested in the Minkowski sum $n \Delta_{d-1}$ of $n$
simplices. Define a \emph{fine mixed cell} of this sum $n
\Delta_{d-1}$ to be a Minkowski sum $B_1 + \cdots + B_n$, where the
$B_i$s are faces of $\Delta_{d-1}$ which lie in independent affine
subspaces, and whose dimensions add up to $d-1$. A \emph{fine mixed
subdivision} of $n\Delta_{d-1}$ is a subdivision\footnote{A
subdivision of a polytope $P$ is a tiling of $P$ with polyhedral
cells whose vertices are vertices of $P$, such that the intersection
of any two cells is a face of both of them.} of $n\Delta_{d-1}$ into
fine mixed cells \cite[Theorem 2.6]{Santos}.

%A \emph{Minkowski cell} of $P$ is a full-dimensional polytope of
%the form $B = B_1 + \cdots + B_k$, where each $B_i$ is a polytope
%with vertices among those of $P_i$. A \emph{mixed subdivision} of
%$P$ is a family of Minkowski cells which cover $P$ and intersect
%properly as Minkowski cells. This means that for any two cells $B
%= \sum B_i$ and $B' = \sum B'_i$, the intersection of $B_i$ and
%$B'_i$ is a face of both. A mixed subdivision is \emph{fine} if,
%for each cell $B = \sum B_i$, the $B_i$s are all simplices and
%$\dim B = \sum \dim B_i$; that is, the $B_i$s lie in independent
%affine subspaces.
%
%%At first sight, the fine mixed subdivisions of a Minkowski sum are
%%quite difficult to work with.
%
%This is not the easiest definition to work with; let us unravel it
%for the case which interests us: the Minkowski sum $n\Delta_{d-1}
%= \Delta_{d-1} + \cdots + \Delta{d-1}$. The cells that we are
%allowed to have in a fine mixed subdivision of $n\Delta_{d-1}$ are
%Minkowski sums $B_1 + \cdots + B_n$ where the $B_i$s are faces of
%$\Delta_{d-1}$ which lie in independent affine subspaces, and
%whose dimensions add up to $d$. Let us call these \emph{fine mixed
%cells}. Also, the cells must intersect properly as Minkowski sums.
%
%Fortunately, Santos showed that any polyhedral subdivision of
%$n\Delta_{d-1}$ into mixed cells can be labelled as a mixed
%subdivision. \cite[Theorem 2.6]{Santos}. Thus instead of checking
%that cells intersect properly as Minkowski sums, we simply have to
%check that they intersect properly; \emph{i.e.}, that the
%intersection of any two cells is a face of both.

\begin{figure}[h]
\centering
\includegraphics[height=5cm]{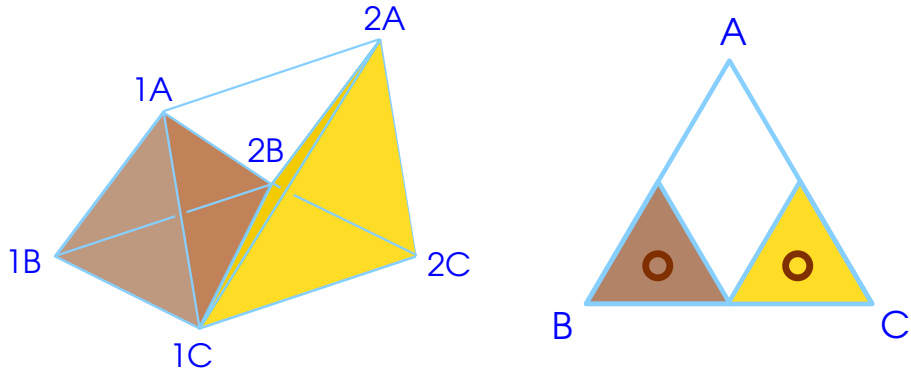}
\caption{The Cayley trick.} \label{fig:cayley.trick}
\end{figure}

Consider the case $d=3$. If the vertices of $\Delta_2$ are
labelled $A,B,$ and $C$, there are two different kinds of fine
mixed cells: a unit triangle like $ABC + A + B + \cdots + A$, and
a unit rhombus like $AB + AC + A + \cdots + C$ (which can face in
three possible directions). Therefore the fine mixed subdivisions
of the triangle $n \Delta_2$ are precisely its lozenge tilings. In
these sums, the summands which are not points determine the shape
of the fine mixed cell, while the summands which are points
determine the position of that cell in the tiling of $n \Delta_2$.
This is illustrated in the right hand side of Figure
\ref{fig:cayley.trick}: a lozenge tiling of $2 \Delta_2$ whose
tiles are $ABC+B$, $AC+AB$, and $C+ABC$.

For $d=4$, if we label the tetrahedron $ABCD$, we have three
different kinds of fine mixed cells: a tetrahedron like $ABCD + A +
\cdots$, a triangular prism like $ABC + AD + A + \cdots$, and a
parallelepiped like $AB + AC + AD + A + \cdots$. These types of
cells correspond to the partitions of $3$ into non-negative
integers: they are $3, 2+1,$ and $1+1+1$, respectively. More
generally, for $\Delta_{d-1}$, there is exactly one type of fine
mixed cell for each partition of $d-1$ into positive integers.

%Throughout this section we will assume some familiarity with
%Minkowski sums of polytopes, mixed subdivisions of Minkowski sums,
%and the Cayley trick, which provides a bijection between fine
%mixed subdivisions of a Minkowski sum $P_1 + \cdots + P_k$ and the
%triangulations of the \emph{Cayley embedding} ${\cal C}(P_1,
%\ldots, P_k)$. For more information on these topics, we refer the
%reader to \cite{Santos}. For a more thorough introduction to
%subdivisions of polytopes, see \cite{deLoera}. {\bf Maybe we
%should say something more about this? A short example for $n
%\Delta_3$ and $\Delta_{n-1} \times \Delta_3$ could be useful.}
%
%We can think of the triangle $T(n)$ as the Minkowski sum
%$\Delta_2+ \cdots + \Delta_2  = n\Delta_2$. The rhombus tilings of
%$n\Delta_2$ are precisely the \emph{fine mixed subdivisions} of
%this Minkowski sum \cite{Santos}.
%
%More generally, we are interested in the fine mixed subdivisions
%of the Minkowski sum $n \Delta_d$. Via the Cayley trick, this is
%equivalent to studying triangulations of the product $\Delta_{n-1}
%\times \Delta_d$. These two points of view have different
%advantages: in the first, we study a simpler, lower-dimensional
%polytope, subdivided into more complicated tiles; in the second,
%the tiles are very simple, but the polytope is more difficult to
%visualize.

In the same way that we identified arrays of triangles with
triangular arrays of dots in Section
\ref{section:rhombus.tilings}, we can identify the array of
possible locations of the simplices in $n \Delta_{d-1}$ with the
array of dots $T_{n,d}$ defined in Section
\ref{section:flags.to.simplex}. A conjectural generalization of
Theorem \ref{thm:tilings}, which we now state, would show that
fine mixed subdivisions of $n \Delta_{d-1}$ are also closely
connected to the matroid $\T_{n,d}$.

\begin{conjecture}\label{conj:subdivisions}
The possible locations of the simplices in a fine mixed
subdivision of $n \Delta_{d-1}$ are precisely the bases of the
matroid $\T_{n,d}$.
\end{conjecture}

In the remainder of this section, we will give a completely
combinatorial description of the fine mixed subdivisions of $n
\Delta_{d-1}$. Then, in Section
\ref{section:triangulations.matroid}, we will use this description
to prove Proposition \ref{prop:subdivisions}, which is the forward
direction of Conjecture \ref{conj:subdivisions}.

We start by recalling the one-to-one correspondence between the fine
mixed subdivisions of $n\Delta_{d-1}$ and the triangulations of
$\Delta_{n-1} \times \Delta_{d-1}$. This equivalent point of view
has the drawback of bringing us to a higher-dimen\-sional picture.
Its advantage is that it simplifies greatly the combinatorics of the
tiles, which are now just
simplices. %This will allow us to write down, in Proposition
%\ref{prop:trees}, a purely combinatorial characterization of the
%pure mixed subdivisions of $n \Delta_{d-1}$.

Let $v_1, \ldots, v_n$ and $w_1, \ldots, w_d$ be the vertices of
$\Delta_{n-1}$ and $\Delta_{d-1}$, so that the vertices of
$\Delta_{n-1} \times \Delta_{d-1}$ are of the form $v_i \times
w_j$. A triangulation $T$ of $\Delta_{n-1} \times \Delta_{d-1}$ is
given by a collection of simplices.  For each simplex $t$ in $T$,
consider the fine mixed cell whose $i$-th summand is $w_aw_b\ldots
w_c$, where $a,b,\ldots,c$ are the indexes $j$ such that $v_i
\times w_j$ is a vertex of $t$. These fine mixed cells constitute
the fine mixed subdivision of $n \Delta_{d-1}$ corresponding to
$T$. (This bijection is only a special case of the more general
Cayley trick, which is discussed in detail in \cite{Santos}.)

For instance, Figure \ref{fig:cayley.trick} (best seen in color) shows a
triangulation of the triangular prism $\Delta_1 \times \Delta_2 = 12
\times ABC$, and the corresponding fine mixed subdivision of $2
\Delta_2$, whose three tiles are $ABC+B, AC+AB$, and $C+ABC$.

\medskip

\begin{figure}[h]
\centering
\includegraphics[height=3.3cm]{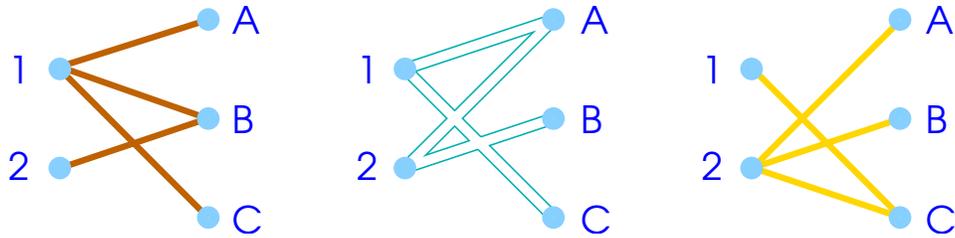}
\caption{The trees corresponding to the triangulation of Figure
\ref{fig:cayley.trick}.} \label{fig:trees}
\end{figure}

Consider the complete bipartite graph $K_{n,d}$ whose vertices are
$v_1, \ldots, v_n$ and $w_1, \ldots, w_d$. Each vertex of
$\Delta_{n-1} \times \Delta_{d-1}$ corresponds to an edge of
$K_{n,d}$. The vertices of each simplex in $\Delta_{n-1} \times
\Delta_{d-1}$ determine a subgraph of $K_{n,d}$. Each
triangulation of $\Delta_{n-1} \times \Delta_{d-1}$ is then
encoded by a collection of subgraphs of $K_{n,d}$. Figure
\ref{fig:trees} shows the three trees that encode the
triangulation of Figure \ref{fig:cayley.trick}.

Our next result is a combinatorial characterization of the
triangulations of $\Delta_{n-1} \times \Delta_{d-1}$.

\newpage

\begin{proposition}\label{prop:trees}
A collection of subgraphs $t_1, \ldots, t_k$ of $K_{n,d}$ encodes
a triangulation of $\Delta_{n-1} \times \Delta_{d-1}$ if and only
if:
\begin{enumerate}
\item Each $t_i$ is a spanning tree.

\item For each $t_i$ and each internal\footnote{An edge of a tree
is \emph{internal} if it is not adjacent to a leaf.} edge $e$ of
$t_i$, there exists an edge $f$ and a tree $t_j$ with $t_j = t_i -
e \cup f$.

\item There do not exist two trees $t_i$ and $t_j$, and a circuit
$C$ of $K_{n,d}$ which alternates between edges of $t_i$ and edges
of $t_j$.
\end{enumerate}
\end{proposition}

\begin{proof}

If $e_1, \ldots, e_n, f_1, \ldots, f_d$ is a basis of $\R^{n+d}$,
then a realization of the polytope $\Delta_{n-1} \times
\Delta_{d-1}$ is given by assigning the vertex $v_i \times w_j$
coordinates $e_i+f_j$. It is then easy to see that the oriented
matroid of affine dependencies of $\Delta_{n-1} \times \Delta_{d-1}$
is the same as the oriented matroid of the graph $K_{n,d}$, with
edges oriented $v_i \rightarrow w_j$ for $1 \leq i \leq n, 1 \leq j
\leq d$. In other words, each minimal affinely dependent set $C$ of
vertices of $\Delta_{n-1} \times \Delta_{d-1}$ corresponds to a
circuit of the graph $K_{n,d}$. Furthermore, the sets $C^+$ and
$C^-$ of vertices which have positive and negative coefficients in
the affine dependence relation of $C$ correspond, respectively, to
the edges that the circuit of $K_{n,d}$ traverses in the forward and
backward direction. Therefore, a set of vertices of $\Delta_{n-1}
\times \Delta_{d-1}$ forms an $(n+d-2)$-dimensional simplex if and
only if it is encoded by a spanning tree of $K_{n,d}$.

%
%each linear equation representing a minimal affine dependence
%between the vertices of $\Delta_{n-1} \times \Delta_{d-1}$
%corresponds to a circuit of the graph $K_{n,d}$. Furthermore, the
%vertices with positive and negative coefficients in the linear
%equation correspond, respectively, to the edges that the circuit
%traverses in the forward and backward direction.

The three conditions in the statement of Proposition
\ref{prop:trees} simply rephrase the following result
\cite[Theorem 2.4.(f)]{SantosAMS}:
\begin{quote}
Suppose we are given a polytope $P$, and a non-empty
collection of simplices whose vertices are vertices of $P$. The
simplices form a triangulation of $P$ if and only if they satisfy
the \emph{pseudo-manifold property}, and \emph{no two simplices
overlap on a circuit}.
\end{quote}

The \emph{pseudo-manifold property} is that, for any simplex
$\sigma$ and any facet $\tau$ of $\sigma$, either $\tau$ is in a
facet of $P$, or there exists another simplex $\sigma'$ with $\tau
\subset \sigma'$. The facets of $\Delta_{n-1} \times \Delta_{d-1}$
are of the form $F \times \Delta_{d-1}$ for a facet $F$ of
$\Delta_{n-1}$ (obtained by deleting one of the vertices of
$\Delta_{n-1}$), or $\Delta_{n-1} \times G$ for a facet $G$ of
$\Delta_{d-1}$ (obtained by deleting one of the vertices of
$\Delta_{d-1}$). Therefore, in the simplex $\sigma$ corresponding
to tree $t$, the facet of $\sigma$ corresponding to $t-e$ is in a
facet of $\Delta_{n-1} \times \Delta_{d-1}$ if and only if $t-e$
has an isolated vertex. So in this case, 2. is equivalent to the
pseudo-manifold property.

Two simplices $\sigma$ and $\sigma'$ are said to \emph{overlap on a
signed circuit} $C=(C^+,C^-)$ of $P$ if $\sigma$ contains $C^+$ and
$\sigma'$ contains $C^-$. The circuits of the polytope $\Delta_{n-1}
\times \Delta_{d-1}$ correspond precisely to the circuits of
$K_{n,d}$, which are alternating in sign. Therefore this condition
is equivalent to 3.
\end{proof}

In light of Proposition \ref{prop:trees}, we will call a
collection of spanning trees satisfying the above properties a
triangulation of $\Delta_{n-1} \times \Delta_{d-1}$.

\medskip

Part of Proposition \ref{prop:trees} is implicit in work of
Kapranov, Postnikov, and Zelevinsky \cite[Section 12]{Postnikov},
and Babson and Billera \cite{Babson}. The latter also gave a
different combinatorial description of the \textbf{regular}
triangulations, which we now describe.

Recall the following geometric method for obtaining subdivisions of
a polytope $P$ in $\R^d$. Assign a height $h(v)$ to each vertex $v$
of $P$, lift the vertex $v$ to the point $(v, h(v))$ in $\R^{d+1}$,
and consider the lower facets of the convex hull of those new points
in $\R^{d+1}$. The projections of those lower facets onto the
hyperplane $x_{d+1}=0$ form a subdivision of $P$. Such a subdivision
is called \emph{regular} or \emph{coherent}.

A regular subdivision of the polytope $\Delta_{n-1}\times
\Delta_{d-1}$ is determined by an assignment of heights to its
vertices. This is equivalent to a weight vector $w$ consisting of
a weight $w_{ij}$ for each edge $ij$ of $K_{n,d}$. Let a
\emph{$w$-weighting} be an assignment $(u,v)$ of vertex weights
$u_1, \ldots, u_n,$ $v_1, \ldots, v_d$ to $K_{n,d}$ such that
$u_i+v_j \geq w_{ij}$ for every edge $ij$ of $K_{n,d}$. Say edge
$ij$ is \emph{$w$-tight} if the equality $u_i+v_j=w_{ij}$ holds;
these edges form the \emph{$w$-tight subgraph} of $(u,v)$. A
subgraph of $K_{n,d}$ is \emph{$w$-tight} if it is the $w$-tight
subgraph of \emph{some} $w$-weighting.

\begin{proposition}\label{prop:regular.trees}\cite{Babson}
Let $w$ be a height vector for $\Delta_{n-1} \times \Delta_{d-1}$
or, equivalently, a weight vector on the edges of $K_{n,d}$. The
regular subdivision corres\-ponding to $w$ consists of the maximal
$w$-tight subgraphs of $K_{n,d}$.
\end{proposition}

Say a weight vector $w$ is \emph{generic} if no circuit of
$K_{n,d}$ has alternating sum of weights equal to $0$. We leave it
to the reader to check, using Proposition
\ref{prop:regular.trees}, that generic weight vectors are
precisely the ones that give rise to regular triangulations.
Hence, if $w$ is generic, the maximal $w$-tight subgraphs of
$K_{n,d}$ are trees, and they satisfy the conditions of Proposition \ref{prop:trees}.
It is an instructive exercise to prove this directly.

\section{Subdivisions of $n\Delta_{d-1}$ and the matroid
$\T_{n,d}$.}\label{section:triangulations.matroid}

Having given a combinatorial characterization of the
triangulations of the polytope $\Delta_{n-1} \times \Delta_{d-1}$
in Proposition \ref{prop:trees}, we are now in a position to prove
the forward direction of Conjecture \ref{conj:subdivisions}, which
relates these triangulations to the matroid $\T_{n,d}$. The
following combinatorial lemma will play an important role in our
proof.

\begin{proposition}\label{prop:posets}
Let $n,d$, and $a_1,\ldots,a_d$ be non-negative integers such that
$a_1 + \cdots +a_d \leq n-1$. Suppose we have a coloring of the
$n(n-1)$ edges of the directed complete graph $K_n$ with $d$
colors, such that each color defines a poset on $[n]$; in other
words,
\begin{enumerate}
\item[(a)] the edges $u \rightarrow v$ and $v \rightarrow u$ have
different colors, and

\item[(b)] if $u \rightarrow v$ and $v \rightarrow w$ have the
same color, then $u \rightarrow w$ has that same color.
\end{enumerate}

Call a vertex $v$ \emph{outgoing} if, for every $i$, there exist
at least $a_i$ vertices $w$ such that $v \rightarrow w$ has color
$i$. Then the number of outgoing vertices is at most
$n-a_1-\cdots-a_d$.
\end{proposition}

\begin{proof}
We have $d$ poset structures on the set $[n]$, and this statement
essentially says that we cannot have ``too many" elements which
are ``very large" in all the posets.

Say there are $x$ outgoing vertices, and let $v$ be one of them.
Let $x_i$ be the number of $i$-colored edges which go from $v$ to
another outgoing vertex, so $x_1 + \ldots + x_d = x-1$.

Consider the $x_1$ outgoing vertices $u_1, \ldots, u_{x_1}$ such
that $v \rightarrow u_j$ is blue. The blue subgraph of $K_n$ is a
poset; so among the $u_j$s we can find a minimal one, say $u_1$,
in the sense that $u_1 \rightarrow u_j$ is not blue for any $j$.
Since $u_1$ is outgoing, there are at least $a_1$ vertices $w$ of
the graph such that $u_1 \rightarrow w$ is blue. This gives us
$a_1$ vertices $w$, other than the $u_i$s, such that $v
\rightarrow w$ is blue. Therefore the blue outdegree of $v$ in
$K_n$ is at least $x_1 + a_1$.

Repeating the same reasoning for the other colors, and summing
over all colors, we obtain:
\begin{eqnarray*}
n-1 &=& \sum_{i=1}^d \left( \textrm{color-}i \textrm{ outdegree of } v \right)\\
& \geq & \sum_{i=1}^d \left(x_i + a_i  \right) \\
& = & x - 1 + \sum_{i=1}^d a_i,
\end{eqnarray*}
which is precisely what we wanted to show.%
%We proceed by induction on $N-a_1-\cdots-a_r$.
%
%The first case to consider is $N-a_1-\cdots-a_r = 1$; suppose that
%we have two outgoing vertices $u$ and $v$. Assume, without loss of
%generality, that the edge $u \rightarrow v$ has color $1$. Then
%the outdegree of $u$ in $G_1$ is greater than the outdegree of
%$v$, and therefore it is at least $a_1+1$. By assumption, the
%outdegree of $u$ in $G_i$ is at least $a_i$ for every other color
%$i$. It follows that the total outdegree of $u$ is at least
%$1+a_1+\cdots+a_r = N$, a contradiction.
%
%Now assume $N-a_1-\cdots-a_r > 1$. We are done if there are no
%outgoing vertices; so assume that $v$ is outgoing. Consider the
%complete graph $K_{N-1}$ obtained by removing vertex $v$, and the
%subgraphs $G'i$ of color $i$ in this graph. Suppose $u \neq v$ is
%outgoing in $K_N$. Then it is outgoing with respect to
%$N-1,r,a_1,...,a_i-1,...,a_r$, but this doesn't do it! FIX THIS!!!
\end{proof}

Notice that the bound of Proposition \ref{prop:posets} is optimal.
To see this, partition $[n]$ into sets $A_1, \ldots, A_d, A$ of
sizes $a_1, \ldots, a_d, n-a_1-\cdots-a_d$, respectively. For each
$i$, let the edges from $A$ to $A_i$ have color $i$. Let the edges
from $A_1$ to $A$ have color $d$, and the edges from the other
$A_i$s to $A$ have color $1$. Pick a linear order for $A$, and let
the edges within $A$ have color $d$ in the increasing order, and
color $1$ in the decreasing order. Pick a linear order for $A_1
\cup \cdots \cup A_d$ where the elements of $A_1$ are the smallest
and the elements of $A_d$ are the largest. Let the edges within
$A_1 \cup \cdots \cup A_d$ have color $d$ in the increasing order,
and color $1$ in the decreasing order. It is easy to check that
this coloring satisfies the required conditions, and it has
exactly $n-a_1-\cdots-a_d$ outgoing vertices.

Also notice that our proof of Proposition \ref{prop:posets}
generalizes almost immediately to the situation where we allow
edges to be colored with more than one color.

\medskip

We have now laid down the necessary groundwork to prove one
direction of Conjecture \ref{conj:subdivisions}.

\begin{proposition}\label{prop:subdivisions}
In any fine mixed subdivision of $n \Delta_{d-1}$,
\begin{enumerate}
\item[(a)] there are exactly $n$ tiles which are simplices, and
\item[(b)] the locations of the $n$ simplices give a basis of the
matroid $\T_{n,d}$.
\end{enumerate}
\end{proposition}

\begin{proof}[Proof of Proposition \ref{prop:subdivisions}]

Let us look back at the way we defined the correspondence between
a triangulation $T$ of $\Delta_{n-1} \times \Delta_{d-1}$ and a
fine mixed subdivision $f(T)$ of $n \Delta_{d-1}$. It is clear
that the simplices $f(t)$ of $f(T)$ arise from those simplices $t$
of $T$ whose vertices are $v_i \times w_1, \ldots, v_i \times w_d$
(for some $i$), and one $v_j \times w_{g(j)}$ for each $j \neq i$.
Furthermore, the location of $f(t)$ in $n \Delta_{d-1}$ is given
by the sum of the $w_{g(j)}$s.

\begin{figure}[h]
\centering
\includegraphics[height=2.5cm]{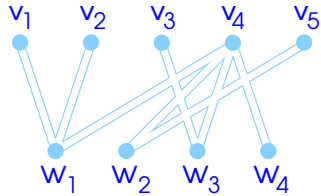}
\caption{A spanning tree of $K_{5,4}$.} \label{fig:tree}
\end{figure}

For instance the spanning tree of $K_{5,4}$ shown in Figure
\ref{fig:tree} gives rise to a simplex in a fine mixed subdivision
of $5\Delta_3 = 5 w_1w_2w_3w_4$ given by the Minkowski sum
$w_1+w_1+w_3+w_1w_2w_3w_4+w_2$. The location of this simplex in
$5\Delta_3$ corresponds to the point $(2,1,1,0)$ of $T_{5,4}$,
because the Minkowski sum above contains two $w_1$ summands, one
$w_2$, and one $w_3$.

In other words, the simplices of the fine mixed subdivision of $n
\Delta_{d-1}$ come from spanning trees $t$ of $K_{n,d}$ for which
one vertex $v_i$ has degree $d$ and the other $v_j$s have degree
$1$. The coordinates of the location of $f(t)$ in $n \Delta_{d-1}$
are simply $(\deg_t w_1-1, \ldots, \deg_t w_d-1)$. Call such a
simplex, and the corresponding tree, \emph{$i$-pure}.
Figure~\ref{fig:tree} shows a 4-pure tree. Also, in the
triangulation of Figures \ref{fig:cayley.trick} and \ref{fig:trees},
there is a $1$-pure tree and a $2$-pure tree, which give simplices
in locations $(0,1,0)$ and $(0,0,1)$ of $2 \Delta_2$, respectively.

\medskip

\noindent \emph{Proof of (a).} We claim that in a triangulation
$T$ of $\Delta_{n-1} \times \Delta_{d-1}$ there is exactly one
$i$-pure simplex for each $i$ with $1 \leq i \leq n$.

First we show there is at least one $i$-pure simplex. If we
restrict the trees of $T$ to the ``claw" subgraph
$K_{\{v_i\},\{w_1,\ldots,w_d\}}$, they should encode a
triangulation of the face $v_i \times (w_1 \ldots w_d)$ of
$\Delta_{n-1} \times \Delta_{d-1}$. This triangulation necessarily
consists of a single simplex, encoded by the claw graph.
Therefore, there must be at least one spanning tree $t$ in $T$
containing this claw.

Now assume that we have two $i$-pure trees $t_1$ and $t_2$. They
must differ somewhere, so assume that $t_1$ contains edge $v_aw_b$
and $t_2$ contains $v_aw_c$. Then we have a circuit $v_aw_bv_iw_c$
of $K_{n,d}$ whose edges alternate between $t_1$ and $t_2$, a
contradiction.

\medskip
\noindent \emph{Proof of (b).} As in the proof of Lemma
\ref{lemma:saturated}, let $T_{a_1,\ldots,a_d}$ be the simplex
consisting of the locations $(x_1, \ldots, x_d)$ in $n
\Delta_{d-1}$ such that $\sum x_i = n$ and $x_i \geq a_i$ for each
$i$. We need to show that $T_{a_1,\ldots,a_d}$, which has a
sidelength of $n-a_1-\cdots-a_d$, contains at most
$n-a_1-\cdots-a_d$ simplices of the fine mixed subdivision.

Somewhat predictably, we will construct a coloring of the directed
complete graph $K_n$ which will allow us to invoke Proposition
\ref{prop:posets}. This coloring will be an economical way of
storing the descriptions of the $n$ pure simplices or,
equivalently, the $n$ pure trees. Let $t_i$ be the $i$-pure tree
in the corresponding triangulation of $\Delta_{n-1} \times
\Delta_{d-1}$. We will color the edge $i \rightarrow j$ in $K_n$
with the color $a$, where $w_a$ is the unique neighbor of vertex
$v_j$ in tree $t_i$. We claim that the $a$-colored subgraph of
$K_n$ is a poset for each color $a$.

First assume that $i \rightarrow j$ and $j \rightarrow i$ have the
same color $a$. Then tree $t_i$ contains edge $v_jw_a$ and tree
$t_j$ contains edge $v_iw_a$. Then, for any $b \neq a$, we have a
circuit $v_iw_bv_jw_a$ of $K_{n,d}$ which alternates between edges
of $t_i$ and $t_j$, a contradiction.

Now assume that $i \rightarrow j$ and $j \rightarrow k$ have color
$a$, but $i \rightarrow k$ has some other color $b$. This means
that $v_jw_a$ and $v_kw_b$ are edges of $t_i$ and $v_jw_a$ is an
edge of $t_j$. But then the circuit $v_jw_av_kw_b$ of $K_{n,d}$
alternates between edges of $t_i$ and $t_j$, a contradiction.

We can now apply Proposition \ref{prop:posets}, and conclude that
there are at most $n-a_1-\ldots-a_d$ outgoing vertices in our
coloring of $K_n$. But observe that the simplex of $n
\Delta_{d-1}$ corresponding to the $i$-pure tree $t_i$ is in
location
%\begin{eqnarray*}
%& & \left((\textrm{degree of }w_1 \textrm { in } t_i) -1 , \ldots,
%\textrm{degree of }w_d \textrm { in } t_i) - 1\right) \\
%& = & (\textrm{color-} 1 \textrm{ outdegree of } i \textrm{ in }
%K_n, \ldots, \textrm{color-} d \textrm{ outdegree of } i \textrm{
%in } K_n).
%\end{eqnarray*}
\[ (\deg_{t_i}(w_1)-1, \ldots, \deg_{t_i}(w_n)-1) = (\textrm{outdeg}_{K_n, \textrm{color } 1}(i), \ldots,
\textrm{outdeg}_{K_n, \textrm{color } d}(i)).
\]

%
%\begin{eqnarray*}
%& & (\deg_{t_i}(w_1)-1, \ldots, \deg_{t_i}(w_n)-1) \\
%& = & (\textrm{outdeg}_{K_n, \textrm{color } 1}(i), \ldots,
%\textrm{outdeg}_{K_n, \textrm{color } d}(i)).
%\end{eqnarray*}

Therefore the simplex of the fine mixed subdivision which
corresponds to $t_i$ is in $T_{a_1,\ldots,a_d}$ if and only if
vertex $i$ is outgoing in our coloring of $K_n$. The desired
result follows.
\end{proof}

For the converse of Conjecture \ref{conj:subdivisions}, we would
need to show that every basis of $\T_{n,d}$ arises from a fine
mixed subdivision of $n \Delta_{d-1}$. In fact, we conjecture a
stronger result, which we state after introducing the necessary
definitions.

Recall the definition of a \emph{regular subdivision} of a polytope
given in Section \ref{section:fine.mixed.subdivisions}. Similarly, a
\emph{regular mixed subdivision} of a Minkowski sum $P_1+ \cdots +
P_n$ in $\R^d$ is obtained by assigning a height $h_i(v)$ to each
vertex $v$ of $P_i$, and projecting the lower facets of the convex
hull of the points in $\R^{d+1}$ of the form $(v_1, h_1(v_1)) +
\cdots + (v_n, h_n(v_n))$, where $v_i$ is a vertex of $P_i$.

%assign a height $h(v)$ to each vertex $v$ of $P$ in such a way
%that $S$ is the projection onto the hyperplane $x_{d+1}=0$ of the
%lower facets of the convex hull of the points $(v, h(v))$ in
%$\R^{d+1}$.

\begin{conjecture}\label{conj:regular.subdivisions}
For any basis $B$ of $\T_{n,d}$, there is a {\bf regular} fine
mixed subdivision of $n \Delta_{d-1}$ whose $n$ simplices are
located at $B$.
\end{conjecture}

The Cayley trick provides us with a bijection between the
triangulations of $\Delta_{n-1} \times \Delta_{d-1}$ and the fine
mixed subdivisions of $n \Delta_{d-1}$. This correspondence also
gives a bijection between regular triangulations of $\Delta_{n-1}
\times \Delta_{d-1}$ and regular fine mixed subdivisions of $n
\Delta_{d-1}$.\cite[Theorem 3.1]{Huber} There is also a
correspondence between the regular triangulations of $\Delta_{n-1}
\times \Delta_{d-1}$ and the combinatorial types of arrangements of
$d$ generic tropical hyperplanes in tropical $(n-1)$-space
\cite{Develin, Santos}.

Just as the combinatorial properties of real hyperplane
arrangements are captured in the theory of oriented matroids,
tropical hyperplane arrangements deserve an accompanying theory of
\emph{tropical oriented matroids}. The discussion of the previous
paragraph suggests that subdivisions of products of two simplices
play the role of tropical oriented matroids, with regular
subdivisions corresponding to realizable tropical oriented
matroids. The multiple appearances of these subdivisions in the
literature are presumably a good indication of the applicability
of tropical oriented matroid theory. Our ability to attack
Conjectures \ref{conj:subdivisions} and
\ref{conj:regular.subdivisions} is one way to measure our progress
on this theory.

\section{Applications to Schubert calculus.}\label{section:schubert}

In this section, we show some of the implications of our work in
the Schubert calculus of the flag manifold. Throughout this
section, we will assume some familiarity with the Schubert
calculus, though we will recall some of the definitions and
conventions that we will use; for more information, see for
example \cite{Fulton, Manivel}. We will also need some of the
results of Eriksson and Linusson \cite{Eriksson1, Eriksson2} and
Billey and Vakil \cite{Billey} on Schubert varieties and
permutation arrays.

The flag manifold $\Fl_n = \Fl_n(\C)$ is a smooth projective
variety which parameterizes the complete flags in $\C^n$. The
\emph{relative position} of any two flags $\E$ and $\F$ in $\Fl_n$
is given by a permutation $w \in S_n$. Let us explain what this
means.

To the permutation $w$, we associate the permutation
matrix\footnote{Notice that this is slightly different from the
usual convention, but it is useful from the point of view of
permutation arrays.} which has a $1$ in the $w(i)$th row of column
$n-i+1$ for $1 \leq i \leq n$. Let $w[i,j]$ be the principal
submatrix with lower right hand corner $(i,j)$, and form an $n
\times n$ table, called a \emph{rank array}, whose entry $(i,j)$
is equal to $\rk \, w[i,j]$. The matrix and rank array associated
to $w=53124 \in S_5$ are shown below.
\[
\left[ \begin{array}{ccccc} 0 & 0 & 1 & 0 & 0 \\
0 & 1 & 0 & 0 & 0 \\
0 & 0 & 0 & 1 & 0 \\
1 & 0 & 0 & 0 & 0 \\
0 & 0 & 0 & 0 & 1
\end{array} \right] \rightarrow
\left[ \begin{array}{ccccc} 0 & 0 & 1 & 1 & 1 \\
0 & 1 & 2 & 2 & 2 \\
0 & 1 & 2 & 3 & 3 \\
1 & 2 & 3 & 4 & 4 \\
1 & 2 & 3 & 4 & 5
\end{array} \right].
\]
Saying that $\E$ and $\F$ are in relative position $w$ means that
the dimensions $\dim(E_i \cap F_j)$ are given precisely by the
rank array of $w$; that is,
\[
\dim(E_i \cap F_j) = \rk \, w[i,j] \qquad \textrm{for all}
\,\,\,\, 1 \leq i,j \leq n.
\]

\medskip

Eriksson and Linusson \cite{Eriksson1, Eriksson2} introduced a
higher-dimensional analog of a permutation matrix, called a
\emph{permutation array}. A permutation array is an array of dots
in the cells of a $d$-dimensional $n \times n \times \cdots \times
n$ box, satisfying some quite restrictive properties. From a
permutation array $P$, via a simple combinatorial rule, one can
construct a \emph{rank array} of integers, also of shape $[n]^d$.
We denote it $\rk \, P$. This definition is motivated by their
result \cite{Eriksson2} that the relative position of $d$ flags
$\E^1, \ldots, \E^d$ in $\Fl_n$ is described by a unique
permutation array $P$, via the equations
\[
\dim\bigl(E^1_{x_1} \cap \cdots \cap E^d_{x_d}\bigr) = \rk \,
P[x_1, \ldots, x_d] \qquad \textrm{for all} \,\,\,\, 1 \leq x_1,
\ldots, x_d \leq n.
\]
This result initiated the study of \emph{permutation array
schemes}, which generalize Schubert varieties in the flag manifold
$\Fl_n$. These schemes are much more subtle than their
counterparts; they can be empty, and are not necessarily
irreducible or even equidimensional. \cite{Billey}

The relative position of $d$ generic flags is described by the
\emph{transversal permutation array}
\[
\Bigl\{(x_1, \ldots, x_d) \in [n]^d \,\, \bigl|\bigr. \,\,
\sum_{i=1}^d x_i = (d-1)n+1 \Bigr\}.
\]
For $\sum_{i=1}^d x_i = (d-1)n+1$, the dot at position $(x_1,
\ldots, x_d)$ represents a one-dimensional intersection $E^1_{x_1}
\cap \cdots \cap E^d_{x_d}$. Naturally, we identify the dots in the
transversal permutation array with the elements of the matroid
$\T_{n,d}$.

\bigskip

Given a fixed flag $\E$, define a \emph{Schubert cell} and
\emph{Schubert variety} to be
\begin{eqnarray*}
X_w^{\circ}(\E) &=& \{\F \, | \, \E \textrm{ and } \F \textrm{
have
relative position } w\} \\
&=& \{\F \, | \, \dim(E_i \cap F_j) = \rk \, w[i,j] \,\,\,
\textrm{for all} \,\,\,\, 1 \leq i,j \leq n.\}, \textrm{ and } \\
X_w(\E) &=& \{\F \, | \, \dim(E_i \cap F_j) \geq \rk \, w[i,j]
\,\,\, \textrm{for all} \,\,\,\, 1 \leq i,j \leq n.\},
\end{eqnarray*}
respectively. The dimension of the Schubert variety $X_w(\E)$ is
$\l(w)$, the number of inversions of $w$.

A \emph{Schubert problem} asks for the number of flags $\F$ whose
relative positions with respect to $d$ given fixed flags $\E^1,
\ldots, \E^d$ are given by the permutations $w^1, \ldots, w^d$.
This question only makes sense when
\[
X = X_{w^1}(\E^1) \cap \cdots \cap X_{w^d}(\E^d)
\]
is $0$-dimensional; that is, when $\l(w^1) + \cdots + \l(w^d) = {n
\choose 2}$. If $\E^1, \ldots, \E^d$ are sufficiently generic, the
intersection $X$ has a fixed number of points $c_{w^1\ldots w^d}$
which only depends on the permutations $w^1, \ldots, w^d$.

This question is a fundamental one for several reasons; the
numbers $c_{w^1\ldots w^d}$ which answer it appear in another
important context. The cycles $[X_w]$ corresponding to the
Schubert varieties form a $\Z$-basis for the cohomology ring of
the flag manifold $\Fl_n$, and the numbers $c_{uvw}$ are the
multiplicative structure constants. (For this reason, if we know
the answer to all Schubert problems with $d=3$, we can easily
obtain them for higher $d$.) The analogous structure constants in
the Grassmannian are the Littlewood-Richardson coefficients, which
are much better understood. For instance, even though the
$c_{uvw}$s are known to be positive integers, it is a long
standing open problem to find a combinatorial interpretation of
them.

\medskip

Billey and Vakil \cite{Billey} showed that the permutation arrays
of Eriksson and Linusson can be used to explicitly intersect
Schubert varieties, and compute the numbers $c_{w^1\ldots w^d}$.

\begin{theorem}(Billey-Vakil, \cite{Billey})\label{thm:billey-vakil}
Suppose that \[ X = X_{w^1}(\E^1) \cap \cdots \cap X_{w^d}(\E^d)
\] is a $0$-dimensional and nonempty intersection, with $\E^1,
\ldots, \E^d$ generic.

\begin{enumerate}
\item There exists a unique permutation array $P \subset
[n]^{d+1}$, easily constructed from $w^1, \ldots, w^d$, such that
\[
\dim\bigl(E^1_{x_1} \cap \cdots \cap E^d_{x_d} \cap F_{x_{d+1}}
\bigr) = \rk \, P[x_1, \ldots, x_d, x_{d+1}],
\]
for all $\F \in X$ and all $1 \leq x_1, \ldots, x_{d+1} \leq n$.

\item These equalities can be expressed as a system of
determinantal equations in terms of the permutation array $P$ and
a vector $v_{a_1, \ldots, a_d}$ in each one-dimensional
intersection $E_{a_1, \ldots, a_d} = E^1_{a_1} \cap \cdots \cap
E^d_{a_d}$. This gives an explicit set of polynomial equations
defining $X$.
\end{enumerate}
\end{theorem}

\medskip

Theorem \ref{thm:billey-vakil} highlights the importance of
studying the line arrangements ${\bf E}_{n,d}$ determined by
intersecting $d$ generic complete flags in $\C^n$. In principle,
if we are able to construct such a line arrangement, we can
compute the structure constants $c_{uvw}$ for any $u,v,w \in S_n$.
(In practice, we still have to solve the system of polynomial
equations, which is not easy for large $n$ or for $d \geq 5$.) Let
us make two observations in this direction.

\subsection{Matroid genericity versus Schubert genericity.}
\label{subsection:genericity}

We have been talking about the line arrangement ${\bf E}_{n,d}$
determined by a generic flag arrangement $\E^1, \ldots, \E^d$ in
$\C^n$. We need to be careful, because we have given two different
meanings to the word \emph{generic}.

In Sections \ref{section:flags.to.simplex}, \ref{section:T_nd} and
\ref{section:T_nd.is.right}, we have shown that, if $\E^1, \ldots,
\E^d$ are sufficiently generic, then the linear dependence
relations in the line arrangement ${\bf E}_{n,d}$ are described by
a fixed matroid $\T_{n,d}$. Let us say that the flags are
\emph{matroid-generic} if this is the case.

Recall that in the Schubert problem described by permutations
$w^1, \ldots, w^d$ with $\sum \l(w^i) = {n \choose 2}$, the
$0$-dimensional intersection
\[
X = X_{w^1}(\E^1) \cap \cdots \cap X_{w^d}(\E^d)
\]
contains a fixed number of points $c_{w^1\ldots w^d}$, provided
that $\E^1, \ldots, \E^d$ are sufficiently generic. Let us say
that $n$ flags in $\C^d$ are \emph{Schubert-generic} if they are
sufficiently generic for any Schubert problem with that given $n$
and $d$.

These notions depend only on the line arrangement ${\bf E}_{n,d}$.
The line arrangement ${\bf E}_{n,d}$ is matroid-generic if its
matroid is $\T_{n,d}$, and it is Schubert-generic if the equations
of Theorem \ref{thm:billey-vakil} give the correct number of
solutions to every Schubert problem.

Our characterization of matroid-generic line arrangements
(\emph{i.e.}, our description of the matroid $\T_{n,d}$) does not
tell us how to construct a Schubert-generic line arrangement.
However, when $d=3$ (which is the interesting case in the Schubert
calculus), the cotransversality of the matroid $\T_{n,3}$ allows
us to present such a line arrangement explicitly.

\begin{proposition} \label{proposition:schubert.generic}
The ${n \choose 2}$ path vectors of Theorem
\ref{thm:represent.cotransversal} are Schubert-generic.
\end{proposition}

\begin{proof}
For each weighting $L$ of the edges of the graph $G_n$ with
complex numbers, like the one shown in Figure \ref{fig:g4}, we can
define the collection $V(L)$ of \emph{path vectors} $v(L)_D =
(v(L)_{D,1}, \ldots, v(L)_{D,n})$ as in Theorem
\ref{thm:represent.cotransversal}: $v(L)_{D,i}$ is the sum of the
weights of all paths from dot $D$ to dot $i$ on the bottom row of
$G_n$.

Consider an arbitrary geometric representation $V$ of $\T_{n,3}$
in $\C^n$. By means of a linear transformation, we can assume that
the vectors assigned to the bottom row are the standard basis
$e_1, \ldots, e_n$, in that order. Say $D$ is any dot in the
triangular array $T_{n,3}$, and $E$ and $F$ are the dots below it.
Since $D,E$ and $F$ are dependent, and $E$ and $F$ are not, we can
write $v_D = ev_E + fv_F$ for some $e,f \in \C$. Write the numbers
$e$ and $f$ on the edges $DE$ and $DF$ of $G_n$. Do this for each
dot $D$, and let $L$ be the resulting weighting of the edges of
$G_n$. Then the collection $V$ is precisely the collection $V(L)$
of path vectors of $L$.

This shows that each matroid-generic line arrangement,
\emph{i.e.}, each geometric representation of $\T_{n,3}$, is given
by the path vectors of a weighting of $G_n$. Among those
matroid-generic line arrangements, the Schubert-generic ones form
an open set, which will include $V(L)$ for any weighting $L$
consisting of algebraically independent weights. This completes
the proof.
\end{proof}

Proposition \ref{proposition:schubert.generic} shows that when we
plug the path vectors $V(L)$ into the polynomial equations of
Theorem \ref{thm:billey-vakil}, and compute the intersection $X$,
we will have $|X| = c_{uvw}$. The advantage of this point of view
is that the equations are now written in terms of combinatorial
objects, without any reference to an initial choice of flags.

\begin{problem}\label{problem:c_uvw}
Interpret combinatorially the $c_{uvw}$ solutions of the above
system of equations, thereby obtaining a combinatorial
interpretation for the structure constants $c_{uvw}$.
\end{problem}

\begin{question}\label{conj:schub.matroid}
Is a Schubert generic flag arrangement always matroid generic?
\end{question}

\begin{question}\label{conj:matroid.schub} Is a
matroid generic flag arrangement always Schubert generic?
\end{question}

We believe that the answer to Question \ref{conj:schub.matroid} is
``yes'' and suspect that the answer to Question
\ref{conj:matroid.schub} is ``no''.

\subsection{A criterion for vanishing Schubert structure
constants.}  \label{subsection:schubert.zero}
%
%Combining the results of \cite{Billey} with our description of
%$\T_{n,d}$, we obtain a very simple and effective criterion for
%proving that many Schubert structure constants are equal to zero.

Consider the Schubert problem
\[
X = X_{w^1}(\E^1) \cap \cdots \cap X_{w^d}(\E^d).
\]
Let $P \in [n]^{d+1}$ be the permutation array which describes the
dimensions $\dim(E^1_{x_1} \cap \cdots \cap E^d_{x_d} \cap
F_{x_{d+1}})$ for any flag $\F \in X$. Let $P_1, \ldots, P_n$ be
the $n$ ``floors" of $P$, corresponding to $F_1, \ldots, F_n$,
respectively. Each one of them is itself a permutation array of
shape $[n]^d$.

Billey and Vakil proposed a simple criterion which is very
efficient in detecting that many Schubert structure constants are
equal to zero.

\begin{proposition}\label{prop:vanishing.billey.vakil}(Billey-Vakil, \cite{Billey})
If $P_n$ is not the transversal permutation array, then
$X=\emptyset$ and $c_{w^1\ldots w^d}=0$.
\end{proposition}

Knowing the structure of the matroid $\T_{n,d}$, we can strengthen
this criterion as follows.

\begin{proposition}\label{prop:vanishing}
Suppose $P_n$ is the transversal permutation array, and identify
it with the set $T_{n,d}$. If, for some $k$, the rank of $P_k \cap
P_n$ in $\T_{n,d}$ is greater than $k$, then $X=\emptyset$ and
$c_{w^1\ldots w^d}=0$.
\end{proposition}

\begin{proof}
Each dot in $P_n$ corresponds to a one-dimensional intersection of
the form $E^1_{x_1} \cap \cdots \cap E^d_{x_d}$. Therefore, each
dot in $P_k \cap P_n$ corresponds to a line that $F_k$ is supposed
to contain if $\F$ is a solution to the Schubert problem. The rank
of $P_k \cap P_n$ is the dimension of the subspace spanned by
those lines; if $\F$ exists, that dimension must be at most $k$.
\end{proof}

Let us see how to apply Proposition \ref{prop:vanishing} in a
couple of examples. Following the algorithm of \cite{Billey}, the
permutations $u=v=w=213$ in $S_3$ give rise to the
four-dimensional permutation array consisting of the dots
$(3,3,1,1)$, $(1,3,3,2)$, $(3,1,3,2)$, $(3,3,1,2)$, $(1,3,3,3)$,
$(2,2,3,3)$, $(2,3,2,3)$, $(3,1,3,3)$, $(3,2,2,3)$, and
$(3,3,1,3)$. We follow \cite{Eriksson2, Vakil} in representing it
as follows:

\[
\tableau{\emt& \emt& \emt \\
\emt& \emt& \emt\\
 \emt& \emt& 1}
\hspace{.3in}
\tableau{\emt& \emt& 3\\
\emt& \emt& \emt\\
3& \emt& 1} \hspace{.3in}
\tableau{\emt& \emt& 3 \\
\emt& 3& 2\\
3& 2& 1}
\]

\medskip

The three boards shown represent the three-dimensional floors
$P_1, P_2,$ and $P_3$ of $P$, form left to right. In each one of
them, a dot in cell $(i,j,k)$ is represented in two dimensions by
a number $k$ in cell $(i,j)$.

It takes some practice to interpret these tables; but once one is
used to them, it is very easy to proceed. Simply notice that $P_2
\cap P_3$ is a set of rank $3$ in the matroid $\T_{3,3}$, while
$P_2$ has rank $2$ as a permutation array, and we are done! We
conclude that $c_{213,213,213}=0$. For $n=3$, this is the only
vanishing $c_{uvw}$ which is not explained by Proposition
\ref{prop:vanishing.billey.vakil}. In this example, the vanishing
of $c_{uvw}$ can also be seen by comparing the leading terms of
the corresponding Schubert polynomials.

\medskip

For a larger example, let $u=2134, v=3142, w=2314$. Notice that
$\l(u) + \l(v) + \l(w) = 6$. The permutation array we obtain is

\[
\tableau{\emt&\emt& \emt& \emt \\
\emt&\emt& \emt& \emt\\
\emt& \emt& \emt& \emt \\
 \emt& 4& \emt& \emt} \hspace{.3in}
\tableau{\emt&\emt& \emt& \emt \\
\emt&\emt& \emt& \emt\\
\emt& \emt& \emt& 4 \\
\emt& 4& \emt& 1} \hspace{.3in}
\tableau{\emt&\emt& \emt& 4 \\
\emt&\emt& \emt& \emt\\
\emt& 4& \emt& 3 \\
4& 3& \emt& 1} \hspace{.3in}
\tableau{\emt&\emt& \emt& 4 \\
\emt&\emt& 4& 3\\
\emt& 4& 3 & 2 \\
4& 3& 2& 1}
\]

\medskip

Here $P_3 \cap P_4$ has rank $4$ in $\T_{4,3}$, implying that
$c_{2134,3142,2314}=0.$

\medskip

Knutson \cite{Knutson}, Lascoux and Schutzenberger \cite{Lascoux}, and
Purbhoo \cite{Purbhoo} have developed other methods for detecting the
vanishing of Schubert structure constants. In comparing these methods
for small values of $n$, we have found Proposition
\ref{prop:vanishing} to be quicker and simpler to observe, but we have
not been able to verify our technique as far as Purbhoo.

Here is an example where our method allows us to ``observe'' a zero
coefficient that Knutson \cite[Fact 2.4]{Knutson} claims does not
follow from his technique of descent cycling.  Let $u=2 3 1 6 4 5$,
$v=2 3 1 6 4 5$, $w=3 2 6 1 5 4$, then the unique permutation array
determined by these three permutations is:

%(setf foo (perm-array '((2 3 1 6 4 5 ) (2 3 1 6 4 5 ) (3 2 6 1 5 4))))
%CL-USER(5): (loop for i from 0 to 5 do (print-3-dim-array-latex (project-array-to-board foo i) 6))
\begin{align}\label{}
&\tableau{
  \emt&\emt&\emt&\emt&\emt&\emt&\\
  \emt&\emt&\emt&\emt&\emt&\emt&\\
  \emt&\emt&\emt&\emt&\emt&\emt&\\
  \emt&\emt&\emt&\emt&\emt&\emt&\\
  \emt&\emt&\emt&\emt&4&\emt&\\
  \emt&\emt&\emt&\emt&\emt&\emt&\\
}
\tableau{
  \emt&\emt&\emt&\emt&\emt&\emt&\\
  \emt&\emt&\emt&\emt&\emt&\emt&\\
  \emt&\emt&\emt&\emt&\emt&\emt&\\
  \emt&\emt&\emt&\emt&5&\emt&\\
  \emt&\emt&\emt&5&4&\emt&\\
  \emt&\emt&\emt&\emt&\emt&\emt&\\
}
\tableau{
  \emt&\emt&\emt&\emt&\emt&\emt&\\
  \emt&\emt&\emt&\emt&\emt&\emt&\\
  \emt&\emt&\emt&\emt&\emt&\emt&\\
  \emt&\emt&\emt&\emt&5&\emt&\\
  \emt&\emt&\emt&5&4&\emt&\\
  \emt&\emt&\emt&\emt&\emt&1&\\
}\\
&
\tableau{
  \emt&\emt&\emt&\emt&\emt&6&\\
  \emt&\emt&\emt&\emt&\emt&\emt&\\
  \emt&\emt&\emt&\emt&\emt&\emt&\\
  \emt&\emt&\emt&\emt&5&\emt&\\
  \emt&\emt&\emt&5&4&\emt&\\
  6&\emt&\emt&\emt&\emt&1&\\
}
\tableau{
  \emt&\emt&\emt&\emt&\emt&6&\\
  \emt&\emt&\emt&\emt&\emt&\emt&\\
  \emt&\emt&\emt&\emt&6&5&\\
  \emt&\emt&\emt&6&5&4&\\
  \emt&\emt&6&5&4&2&\\
  6&\emt&5&4&2&1&\\
}
\tableau{
  \emt&\emt&\emt&\emt&\emt&\cinum{6}&\\
  \emt&\emt&\emt&\emt&6&5&\\
  \emt&\emt&\emt&6&\cinum{5}&4&\\
  \emt&\emt&6&\cinum{5}&4&3&\\
  \emt&6&\cinum{5}&4&3&2&\\
  \cinum{6}&5&4&3&2&\cinum{1}&\\
}
\end{align}
By Theorems~\ref{thm:matroid} and~\ref{thm:genericflags}, an
independent set on the rank 6 board is determined by the circled
points
\[
\{(1,6,6), (3,5,5), (4,4,5), (5,3,5), (6,1,6), (6,6,1) \}.
\]
In the rank 4 board, we have the points $(1,6,6), (6,1,6), (6,6,1)$
from this basis along with $(4,5,5), (5,4,5),(5,5,4)$ which span a
two dimensional space in the span of $(3,5,5), (4,4,5), (5,3,5)$.
Therefore, the 4 dimensional board cannot be satisfied by vectors in
a space of dimension less than 5.  Hence $c_{uvw}=0$.

Proposition \ref{prop:vanishing} is only the very first observation
that we can make from our understanding of the structure of
$\T_{n,d}$. Our argument can be easily fine-tuned to explain all
vanishing Schubert structure constants with $n \leq 5$. A systematic
way of doing this in general would be very desirable.

\section{Future directions.}\label{sec:future}

We invite our readers to pursue some further directions of study
suggested by the results in this paper. Here they are, in order of
appearance.

\begin{itemize}

\begin{figure}[h]
\centering
\includegraphics[height=3cm]{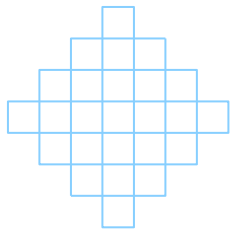}
\caption{A fool's diamond.} \label{fig:aztec}
\end{figure}

\item Theorem \ref{thm:tilings} generalizes to rhombus tilings of
any region in the triangular lattice, or domino tilings of any
region in the square lattice. If $R$ is a region with more upward
than downward triangles (or more black than white squares), let
$\B$ be the sets of holes such that the remaining holey region $R$
has a rhombus tiling (or a domino tiling). Then $\B$ is the set of
bases of a matroid $M_R$. Are there any other regions $R$ for
which the matroid $M_R$ has a nice geometric interpretation? A
good candidate, suggested by Jim Propp, is what he calls the
\emph{fool's diamond} \cite{Propp}, shown in Figure
\ref{fig:aztec}.

\item Subdivisions of $\Delta_{n-1} \times \Delta_{d-1}$, or
equivalently tropical oriented matroids, appear in many different
contexts. A detailed investigation of these objects promises to
become a useful tool. Aside from their intrinsic interest,
Conjectures \ref{conj:subdivisions} and
\ref{conj:regular.subdivisions} should help us develop this
tropical oriented matroid theory.

%Conjectures \ref{conj:subdivisions} and
%\ref{conj:regular.subdivisions} are still open. Aside from their
%intrinsic interest, they are a good starting point in our search
%for a better understanding of the combinatorics of triangulations
%of products of simplices, and of ``tropical oriented matroids".

\item We still do not have a solid understanding of the
relationship between two of the main subjects of our paper: the
geometry of $d$ flags in $\C^n$ and the triangulations of
$\Delta_{n-1} \times \Delta_{d-1}$. We have shown that some
aspects of the geometric information of the flags (the
combinatorics of the line arrangement they determine, and the
vanishing of many Schubert structure constants) are described in a
small set of tiles of the triangulations (the $n$ ``pure" tiles).
Can we use the complete triangulations and fine mixed subdivisions
to understand more subtle geometric questions about flags? Does
the geometry of flags tell us something new about triangulations
of products of simplices, and their multiple appearances in
tropical geometry, optimization, and other subjects?

\item In particular, do the triangulations of $\Delta_{n-1} \times
\Delta_{d-1}$ play a role in the Schubert calculus of the flag
manifold $\Fl_n$? Is this point of view related to Knutson, Tao,
and Woodward's use of puzzles \cite{Knutson1, Knutson2} in the
Grassmannian Schubert calculus? Readers familiar with puzzles may
have noticed the similarities and the differences between them and
lozenge tilings of triangles.

\item Problem \ref{problem:c_uvw} is a promising way of attacking
the long-standing open problem of interpreting $c_{uvw}$
combinatorially.

\item Questions~\ref{conj:schub.matroid} and \ref{conj:matroid.schub} remain open.  Is a Schubert
generic flag arrangement always matroid generic?  Is a matroid generic
flag arrangement always Schubert generic?

\item Proposition \ref{prop:vanishing} is just the first
consequence of the matroid $\T_{n,d}$ on the vanishing of the
Schubert structure constants. This argument can be extended in many
ways to explain why other $c_{uvw}$s are equal to $0$. A systematic
way of doing this would be desirable, and seems within reach at
least for $n \leq 7$.

\end{itemize}

\section{Acknowledgments}

We would like to thank Laci Lovasz, Jim Propp, Kevin Purbhoo, Ravi
Vakil, David Wilson, and Andrei Zelevinsky for very helpful
discussions.

\small{

}

\begin{flushleft}
\bigskip
\bigskip
\footnotesize{\textsc{Department of Mathematics \\San Francisco
State University \\ San Francisco, CA, USA} \\
\texttt{federico@math.sfsu.edu}

\bigskip

\textsc{Department of Mathematics \\University of Washington \\
Seattle, WA, USA} \\
\texttt{billey@math.washington.edu} }

\end{flushleft}

\end{document}